\documentclass[12pt]{article}
\usepackage{amsthm,amsfonts,amssymb,amscd}

\begin{document}

\begin{center}
\large \bf Birational geometry \\
of Fano hypersurfaces of index two
\end{center}\vspace{0.5cm}

\centerline{Aleksandr V. Pukhlikov}\vspace{0.5cm}

\parshape=1
3cm 10cm \noindent {\small \quad\quad\quad \quad\quad\quad\quad
\quad\quad\quad {\bf }\newline We prove that every non-trivial
structure of a rationally connected fibre space on a generic (in
the sense of Zariski topology) hypersurface $V$ of degree $M$ in
the $(M+1)$-dimensi\-onal projective space for $M\geq 14$ is given
by a pencil of hyperplane sections. In particular, the variety $V$
is non-rational and its group of birational self-maps coincides
with the group of biregular automor\-phisms and for that reason is
trivial. The proof is based on the techniques of the method of
maximal singularities and inversion of adjunction.

Bibliography: 22 titles.} \vspace{1cm}

\section*{Introduction}

{\bf 0.1. Statement of the main result.} Fix an integer $M\geq 4$.
Denote by the symbol ${\mathbb P}$ the complex projective space
${\mathbb P}^{M+1}$. Let $V=V_M\subset{\mathbb P}$ be a
non-singular hypersurface of degree $M$. Obviously, $V$ is a Fano
variety of index two:
$$
\mathop{\rm Pic}V={\mathbb Z}H,\quad K_V=-2H,
$$
where $H$ is the class of a hyperplane section. On the variety $V$
there are the following structures of a non-trivial rationally
connected fibre space: let $P\subset{\mathbb P}$ be an arbitrary
subspace of codimension two, $\alpha_P\colon{\mathbb
P}\dashrightarrow {\mathbb P}^1$ the corresponding linear
projection, then its restriction
$$
\pi_P=\alpha_P|_V\colon V\dashrightarrow {\mathbb P}^1
$$
fibres $V$ into Fano hypersurfaces of index one and for that
reason defines on $V$ a structure of rationally connected fibre
space. Recall \cite{Pukh13a}, that a (non-trivial) rationally
connected fibre space is a surjective morphism $\lambda\colon Y\to
S$ of projective varieties, where $\mathop{\rm dim}S\geq 1$ and
the variety $S$ and the fibre of general position
$\lambda^{-1}(s)$, $s\in S$, are rationally connected (and the
variety $Y$ itself is automatically rationally connected by the
theorem of Graber, Harris and Starr \cite{GHS}).\vspace{0.1cm}

Here is the main result of the present paper.\vspace{0.1cm}

{\bf Theorem 1.} {\it Assume that $M\geq 14$ and the hypersurface
$V$ is sufficiently general (in the sense of Zariski topology on
the space of coefficients of homogeneous polynomials of degree $M$
on ${\mathbb P}$). Let $\chi\colon V\dashrightarrow Y$ be a
birational map onto the total space of a rationally connected
fibre space $\lambda\colon Y\to S$. Then $S={\mathbb P}^1$ and for
some isomorphism $\beta\colon{\mathbb P}^1\to S$ and some subspace
$P\subset{\mathbb P}$ of codimension two we have
$$
\lambda\circ\chi=\beta\circ\pi_P,
$$
that is, the following diagram commutes:}
$$
\begin{array}{ccccc}
& V & \stackrel{\chi}{\dashrightarrow} & Y &\\
\pi_P & \downarrow & & \downarrow & \lambda\\
& {\mathbb P}^1 & \stackrel{\beta}{\to} & S. &\\
\end{array}
$$
\vspace{0.1cm}

{\bf Corollary 1.} {\it For a generic hypersurface $V$ of
dimension $\mathop{\rm dim} V\geq 14$ the following claims
hold.}\vspace{0.1cm}

(i) {\it On the variety $V$ there are no structures of a
rationally connected fibre space with the base of dimension $\geq
2$. In particular, on $V$ there are no structures of a conic
bundle and del Pezzo fibration, and the variety $V$ itself is
non-rational.\vspace{0.1cm}

{\rm (ii)} Assume that there is a birational map $\chi\colon
V\dashrightarrow Y$, where $Y$ is a Fano variety of index $r\geq
2$ with factorial terminal singularities, such that $\mathop{\rm
Pic}Y={\mathbb Z}H_Y$, where $K_Y=-rH_Y$, and the linear system
$|H_Y|$ is non-empty and free. Then $r=2$ and the map $\chi$ is a
biregular isomorphism.\vspace{0.1cm}

{\rm (iii)} The group of birational self-maps of the variety $V$
coincides with the group of biregular automorphisms:
$$
\mathop{\rm Bir}V=\mathop{\rm Aut}V
$$
and for that reason is trivial.} \vspace{0.1cm}

{\bf Proof of the corollary.} The claims (i-iii) follow from
Theorem 1 in an obvious way. Q.E.D.\vspace{0.1cm}

{\bf Conjecture 1.} {\it Assume that $V_d\subset {\mathbb P}$ is a
smooth hypersurface of degree $d\leq M$, where $d\geq [(M+5)/2]$
(in that case $V_d$ is a Fano variety of index $r=M+2-d$). Let
$\chi\colon V\dashrightarrow Y$ be a birational map onto the total
space of a rationally connected fibre space $\lambda\colon Y\to
S$. Then $\mathop{\rm dim} S\leq r-1$ and if $\mathop{\rm dim} S=
r-1$, then there is a linear subspace $P\subset{\mathbb P}$ of
codimension $r$ and a birational map $\beta\colon{\mathbb
P}^{r-1}\dashrightarrow S$ such that
$\lambda\circ\chi=\beta\circ\pi_P$, that is, the following diagram
commutes}
$$
\begin{array}{ccccc}
& V_d & \stackrel{\chi}{\dashrightarrow} & Y &\\
\pi_P & \downarrow & & \downarrow & \lambda\\
& {\mathbb P}^{r-1} & \stackrel{\beta}{\dashrightarrow} & S. &\\
\end{array}
$$
\vspace{0.1cm}

{\bf Remark 1.} For $d\leq M-1$ (that is, for $r\geq 3$) one can
certainly not expect that {\it all} structures of a rationally
connected fibre space (or of a Fano-Mori fibre space) are linear
projections. Already for a hypersurface of index 3 every pencil of
quadrics defines a rational map onto ${\mathbb P}^1$, the fibre of
which is a complete intersection of the type $2\cdot (M-1)$ in
${\mathbb P}^{M+1}$, that is, a Fano variety of index
one.\vspace{0.1cm}

The purpose of the present paper is to prove Theorem 1. As usual,
its claim will be derived from a lot more technical and less
visual description of maximal singularities of mobile linear
systems on $V$. However, before explaining the structure of the
proof of Theorem 1, let us give a precise meaning to the
assumption of the hypersurface $V$ being generic in the sense of
Zariski topology.\vspace{0.3cm}

%%%%%%%%%%%%%%%%%%%%%%%%%%%%%%%%%%%%%%%%%%%%%%%%%%%%%%%%%%%%%%
%%%%%%%%%%%%%%%%%%%%%%%%%%%%%%%%%%%%%%%%%%%%%%%%%%%%%%%%%%%%%%
%%%%%%%%%%%%%%%%%%%%%%%%%  subsection 0.2

{\bf 0.2. The regularity conditions.} Let
$$
{\cal F}={\mathbb P}(H^0({\mathbb P},{\cal O}_{\mathbb P}(M)))
$$
be the space parametrizing hypersurfaces of degree $M$ in
${\mathbb P}$. The local {\it regularity conditions}, given below,
define an open subset ${\cal F}_{\rm reg}\subset {\cal F}$. A
separate (but not difficult) problem is to show that for $M\geq
14$ the set ${\cal F}_{\rm reg}$ is non-empty.\vspace{0.1cm}

Let $o\in{\mathbb P}$ be an arbitrary point,
$(z_1,\dots,z_{M+1})=(z_*)$ a system of affine coordinates with
the origin at the point $o$ and $V\ni o$ a non-singular
hypersurface of degree $M$. It is given by an equation $f=0$,
where
$$
f=q_1+q_2+\dots +q_M
$$
is a non-homogeneous polynomial in the variables $z_*$, $q_i$ is
its homogeneous component of degree $i$. Let $\Pi\subset{\mathbb
C}^{M+1}$ be an arbitrary linear subspace of codimension
$c\in\{0,1,2,3\}$, on which $q_1$ does not vanish identically,
that is, $\Pi\not\subset T_oV$. We will need the following
regularity conditions.\vspace{0.1cm}

(R1) For any subspace $\Pi$ the sequence of polynomials
$$
q_1|_{\Pi},\, q_2|_{\Pi},\,\dots \,q_{M-c}|_{\Pi}
$$
is regular in ${\cal O}_{o,\Pi}$, that is the system of equations
$$
q_1|_{\Pi}= q_2|_{\Pi}=\dots =q_{M-c}|_{\Pi}=0
$$
determines a finite set of lines.\vspace{0.1cm}

(R2) The rank of the quadratic form
$$
q_2|_{\{q_1=0\}}
$$
is at least $M-[\sqrt{M}]$.\vspace{0.1cm}

(R3) The restriction of the equation $q_3=0$ onto the quadric
hypersurface $\{q_2|_{\Lambda}=0\}$, where $\Lambda$ is an
arbitrary linear subspace of codimension two in the tangent
hyperplane, defines an irreducible reduced closed
set.\vspace{0.1cm}

The last (forth) regularity condition is a global
one.\vspace{0.1cm}

(R4) The intersection of the hypersurface $V$ with an arbitrary
linear subspace $P\subset{\mathbb P}$ of codimension two has at
most isolated quadratic singularities.\vspace{0.1cm}

The following claim is true.\vspace{0.1cm}

{\bf Theorem 2.} {\it For $M\geq 14$ there exists a non-empty
Zariski open subset ${\cal F}_{\rm reg}\subset {\cal F}$, such
that every hypersurface $V\in{\cal F}_{\rm reg}$ is non-singular
and satisfies the conditions (R1-R3) at every point, and also the
condition (R4).}\vspace{0.1cm}

For the {\bf proof} of Theorem 2 see Subsection 1.6.\vspace{0.3cm}

%%%%%%%%%%%%%%%%%%%%%%%%%%%%%%%%%%%%%%%%%%%%%%%%%%%%%%%%%%%%%%
%%%%%%%%%%%%%%%%%%%%%%%%%%%%%%%%%%%%%%%%%%%%%%%%%%%%%%%%%%%%%%
%%%%%%%%%%%%%%%%%%%%%%%%%  subsection 0.3

{\bf 0.3. Plan of the proof of Theorem 1.} For an arbitrary
subspace $P\subset{\mathbb P}$ of codimension two denote by the
symbol $V_P$ the blow up of $V$ along the subvariety $V\cap P$.
For a mobile linear system $\Sigma$ on $V$ its strict transform on
$V_P$ denote by the symbol $\Sigma_P$. Considering instead of
$\Sigma$ its symmetric square, we can alwasy assume that
$\Sigma\subset |2nH|$. Recall that $c_{\rm virt}(\Sigma)$ is the
virtual threshold of canonical adjunction \cite[Sec.
2.1]{Pukh13a}. Theorem 1 is an easy corollary from the technical
fact formulated below.\vspace{0.1cm}

{\bf Theorem 3.} {\it Assume that $M\geq 14$ and $V\in {\cal
F}_{\rm reg}$. If the mobile system $\Sigma\subset |2nH|$
satisfies the inequality
\begin{equation}\label{24.09.2}
c_{\rm virt}(\Sigma)<n,
\end{equation}
then there exists a unique linear subspace $P\subset {\mathbb P}$
of codimension two, such that the subvariety $B=P\cap V$ satisfies
the inequality
\begin{equation}\label{24.09.1}
\mathop{\rm mult}\nolimits_B \Sigma>n,
\end{equation}
whereas for the strict transform $\Sigma_P$ the following equality
holds:}
$$
c_{\rm virt}(\Sigma)=c_{\rm virt}(\Sigma_P)=c(V_P,\Sigma_P).
$$

The system $\Sigma$ and the integer $n\geq 1$ are fixed throughout
the paper. In its turn, Theorem 3 will be derived from the
following two key facts.\vspace{0.1cm}

{\bf Theorem 4.} {\it Assume that for some subvariety $B\subset V$
of codimension two the inequality (\ref{24.09.1}) holds. Then
$B=P\cap V$, where $P\subset{\mathbb P}$ is a linear subspace of
codimension two.}\vspace{0.1cm}

{\bf Theorem 5.} {\it Assume that the inequality (\ref{24.09.2})
holds. Then for some irreducible subvariety $B$ of codimension two
the inequality (\ref{24.09.1}) holds.}\vspace{0.1cm}

Theorems 4 and 5 are given in the order in which they are shown.
Theorem 5 (the exclusion of the infinitely near case) is the most
difficult to prove. Further work is organized as
follows.\vspace{0.1cm}

In Sec. 1, assuming Theorem 3, we show Theorem 1, and after that,
obtain Theorem 3, assuming Theorems 4 and 5. In Sec. 2 we show
Theorem 4. In Sec. 3-5 we prove Theorem 5.\vspace{0.3cm}

%%%%%%%%%%%%%%%%%%%%%%%%%%%%%%%%%%%%%%%%%%%%%%%%%%%%%%%%%%%%%%%%%%%%
%%%%%%%%%%%%%%%%%%%%%%%%%%%%%%%%%%%%%%%%%%%%%%%%%%%%%%%%%%%%%%%%%%%%
%%%%%%%%%%%%%%%%%   subsection 0.4

{\bf 0.4. Historical remarks.} The result, completely similar to
Theorem 1, has been shown for Fano double spaces of index two in
\cite{Pukh10}, see also Chapter 8 in \cite{Pukh13a}. Prior to the
paper \cite{Pukh10}, the only result giving a complete description
of the structures of a rationally connected fibre space on a Fano
variety of index two, was Grinenko's theorem \cite{Grin03b,Grin04}
on the Veronese double cone, a very special Fano
three-fold.\vspace{0.1cm}

A series of important results on birational geometry of Fano
varieties of index two and higher was obtained by other methods:
by the transcendent method of Clemens and Griffiths \cite{CG} and
its subsequent generalizations (see \cite{Clem82}), and also by
means of Koll\' ar's technique \cite{Kol95a,Kol00}. For the
details, see the introduction to the paper \cite{Pukh10}, where,
in particular, the dramatic story of studying the birational
geometry of the Veronese double cone and (not completed to this
day) studying of the double space of index two is
described.\vspace{0.1cm}

Note that the problem of description of the birational type of
Fano varieties of index higher than one was discussed already in
the classical paper \cite{IM}; Fano himself also worked on the
problem (for the cubic three-fold $V_3\subset{\mathbb P}^4$)
\cite{Fano47}.

%%%%%%%%%%%%%%%%%%%%%%%%%%%%%%%%%%%%%%%%%%%%%%%%%%%%%%%%%%%%%%%%%%%
%%%%%%%%%%%%%%%%%%%%%%%%%%%%%%%%%%%%%%%%%%%%%%%%%%%%%%%%%%%%%%%%%%%
%%%%%%%%%%%%%%%%%%%%%%%%%%%%%%%%%%%%%%%%%%%%%%%%%%%%%%%%%%%%%%%%%%%
%%%%%%%%%%%%%%%%%%%%%%%%%%%%%%%%%%%%%%%%%%%%%%%%%%%%%%%%%%%%%%%%%%%
%%%%%%%%%%%%%%%%%%%%%%%%%%%%   SECTION 1

\section{Pencils of hyperplane sections}

In this section, we prove Theorem 1, assuming the claim of Theorem
3. After Theorem 3 is obtained from Theorems 4 and 5. Finally, we
discuss the (routine) proof of Theorem 2.\vspace{0.3cm}

{\bf 1.1. Fano fibre spaces over ${\mathbb P}^1$.} Let us prove
Theorem 1. Let $\Sigma\subset |2nH|$ be the strict transform on
$V$ of a free linear system on $W$, which is the $\lambda$-pull
back of a very ample linear system on the base $S$. Then the
inequality (\ref{24.09.2}) holds, because $c_{\rm
virt}(\Sigma)=0$. The system $\Sigma\subset |2nH|$ is now fixed.
Assuming the claim of Theorem 3, consider the subspace
$P\subset{\mathbb P}$ of codimension two, such that for $B=P\cap
V$ the inequality (\ref{24.09.1}) holds. Let $\varphi\colon V^+\to
V$ be the blow up of the subvariety $B$ and
$E_B=\varphi^{-1}(B)\subset V^+$ the exceptional
divisor.\vspace{0.1cm}

{\bf Lemma 1.1.} (i) {\it The variety $V^+$ is factorial and has
at most finitely many isolated double points (not necessarily
non-degenerate).\vspace{0.1cm}

{\rm (ii)} The linear projection $\pi_{\mathbb P}\colon{\mathbb
P}\dashrightarrow{\mathbb P}^1$ from the subspace $P$ generates
the regular projection
$$
\pi=\pi_{\mathbb P}\circ\varphi\colon V^+\to{\mathbb P}^1,
$$
the generic fibre of which $F_t=\pi^{-1}(t)$, $t\in{\mathbb P}^1$
is a non-singular Fano variety of index one, and a finite number
of singular fibres have isolated double points.\vspace{0.1cm}

{\rm (iii)} The following equalities hold:
$$
\mathop{\rm Pic}V^+={\mathbb Z}H\oplus{\mathbb Z}E_B={\mathbb
Z}K^+\oplus{\mathbb Z}F,
$$
where $H=\varphi^*H$ for simplicity of notations, $K^+=K_{V^+}$ is
the canonical class of the variety $V^+$, $F$ is the class of the
fibre of the projection $\pi$, where}
$$
K^+=-2H+E,\,\,F=H-E.
$$
\vspace{0.1cm}

{\bf Proof.} These claims are obvious by the regularity conditions
and the well known factoriality of an isolated hypersurface
singularity in the dimension 4 and higher, see \cite{CL}. Q.E.D.
for the lemma.\vspace{0.1cm}

Let $\Sigma^+$ be the strict transform of the system $\Sigma$ on
$V^+$.\vspace{0.1cm}

{\bf Lemma 1.2.} {\it The linear system $\Sigma^+$ is composed
from the pencil} $|F|$: $\Sigma^+\subset |2nH|$.\vspace{0.1cm}

{\bf Proof.} For some $m\in {\mathbb Z}_+$ and $l\in {\mathbb Z}$
we have:
$$
\Sigma^+\subset |-mK^++lF|,
$$
where $m=2n-\mathop{\rm mult}\nolimits_B \Sigma$ and
$l=2(\mathop{\rm mult}\nolimits_B \Sigma -n)\geq 2$. Thus the
threshold of canonical adjunction is
$$
c(\Sigma^+,V^+)=m.
$$
By Theorem 3, $c(\Sigma^+,V^+)=c_{\rm virt}(\Sigma)=0$, so that
$m=0$ and $l=2n$, as we claimed. Q.E.D. for the
lemma.\vspace{0.1cm}

Therefore, the mobile linear system $\Sigma$ is composed from the
pencil of hyperplane sections, containing $B$, which completes the
proof of Theorem 1.\vspace{0.3cm}

%%%%%%%%%%%%%%%%%%%%%%%%%%%%%%%%%%%%%%%%%%%%%%%%%%%%
%%%%%%%%%%%%%%%%%%%%%%%%%%%%%%%%%%%%%%%%%%%%%%%%%%%%
%%%%%%%%%%%%%%%%%   subsection 1.2

{\bf 1.2. Mobile systems on the variety $V$.} Assume the claims of
Theorems 4 and 5. Let us prove Theorem 3. In the notations of
Subsection 1.1 we have to show that for the mobile linear system
$$
\Sigma^+\subset |-mK^++lF|
$$
with $l\in{\mathbb Z}_+$ the equality
$$
c_{\rm virt}(\Sigma)^+=c(\Sigma^+,V^+)=m
$$
holds. (This is precisely the claim of Theorem 3.) Assume the
converse:
$$
c_{\rm virt}(\Sigma)^+ < m,
$$
then the pair $(V^+,\frac{1}{m}\Sigma^+$ is not canonical, that
is, the linear system $\Sigma^+$ has a maximal singularity. Since
every fibre of the fibre space $\pi\colon V^+\to {\mathbb P}^1$ is
a factorial birationally rigid variety, the centre of every
maximal singularity is contained in some fibre $F_t=\pi^{-1}(t)$.
Restricting the linear system $\Sigma^+$ onto such a fibre
$F=F_t$, we obtain an effective divisor $D\in |-mK_F|$, such that
the pair
\begin{equation}\label{25.09.1}
(F,\frac{1}{m}D)
\end{equation}
is not canonical (in fact, not log canonical, but we do not use
that). For any curve $C\subset F$, $C\cap \mathop{\rm Sing}
F=\emptyset$, it is known, see \cite[Chapter 2]{Pukh13a}, that
$\mathop{\rm mult}\nolimits_C D\leq m$, which implies that the
centre of every non canonical singularity of the pair
(\ref{25.09.1}) is either a point, or a curve, passing through a
singularity of $F$. Furthermore, it is well known \cite[Chapter
7]{Pukh13a},  that a smooth point can not be the centre of a non
canonical singularity, and the proof of that fact excludes also
the case when the centre is a curve (since $F$ has only isolated
singularities). Therefore, we may assume that the centre of a
maximal (non canonical) singularity of the pair (\ref{25.09.1}) is
a singular point $o$.\vspace{0.1cm}

At this moment, and up to the end of this section, it is
convenient to slightly change the notations. We denote the variety
$F$ by the symbol $W$. It is a hypersurface of degree $M$ in
${\mathbb P}^M$ with an isolated quadratic point $o\in W$. On $W$
there is an effective divisor $D\sim mH$, where $H$ is the class
of a hyperplane section of $W$, such that the pair
$(W,\frac{1}{m}D)$ is not canonical at the point $o$. We have to
show that this is impossible, that is, to obtain a contradiction.
We do it in several steps, modifying the proof in
\cite{Pukh09b}.\vspace{0.3cm}

%%%%%%%%%%%%%%%%%%%%%%%%%%%%%%%%%%%%%%%%%%%%%%%%%%%%%%%%%%%
%%%%%%%%%%%%%%%%%%%%%%%%%%%%%%%%%%%%%%%%%%%%%%%%%%%%%%%%%%%
%%%%%%%%%%%%%%%%%%%%   subsection 1.3

{\bf 1.3. Step 1: effective divisors on quadrics.} Let
$Q\subset{\mathbb P}^{M-1}$ be an irreducible quadric hypersurface
of rank $\geq 5$, $H_Q\in\mathop{\rm Pic}Q={\mathbb Z}H_Q$ the
class of a hyperplane section and $B\subset Q$ an irreducible
subvariety, which is not contained entirely in $\mathop{\rm
Sing}Q$.\vspace{0.1cm}

{\bf Definition 1.1.} We say that the effective divisor $D$ on $Q$
satisfies the condition $H(m)$ with respect to $B$, where $m\geq
1$ is a fixed integer, if for any point of general position $p\in
B$ (in particular, $p\not\in\mathop{\rm Sing}Q$) there exists a
hyperplane $F(p)\subset E_p$ in the exceptional divisor
$E(p)=\varphi^{-1}_p(p)$ of the blow up $\varphi_p\colon Q_p\to Q$
of the point $p$, such that the inequality
$$
\mathop{\rm mult}\nolimits_pD+\mathop{\rm
mult}\nolimits_{F(p)}\widetilde{D}>2m
$$
holds, where $\widetilde{D}\subset Q_P$ is the strict transform of
the divisor $D$.\vspace{0.1cm}

Note that the divisor $D$ is not assumed to be irreducible, and
the integer $m$ does not depend on the point $p$. It is assumed
that the hyperplane $F(p)$ depends algebraically on the point $p$.
Let $l\geq 1$ be the degree of the hypersurface in ${\mathbb
P}^{M-1}$, which cuts out $D$ on $Q$, that is, $D\sim
lH_Q$.\vspace{0.1cm}

Now, repeating the proof of Proposition 2.1 in \cite{Pukh09b} word
for word, we obtain\vspace{0.1cm}

{\bf Proposition 1.1.} {\it Assume that the inequality
$$
\mathop{\rm dim}B+\mathop{\rm rk}E\geq M+3
$$
holds. Assume, moreover, that an effective divisor $D$ satisfies
the condition $H(m)$ with respect to $B$. Then the following
alternative takes place:\vspace{0.1cm}

(1) either the inequality $l>2m$ holds (and we say that this is
the simple case),\vspace{0.1cm}

(2) or there is a hyperplane section $Z\subset Q$, which contains
entirely the subvariety $B$, such that for a point of general
position $p\in B$ in the notations above
$$
F(p)=\widetilde{Z}\cap E_p,
$$
where $\widetilde{Z}\subset Q_p$ is the strict transform of $Z$ on
$Q_p$, and moreover, $Z$ is contained in the divisor $D$ with the
multiplicity
$$
a>2m-l
$$
(in other words, $D=aZ+D^*$, where the effective divisor $D^*$
does not contain $Z$ as a component; this case we say to be the
hard one).}\vspace{0.1cm}

{\bf Remark 1.1.} If the quadric $E$ is non-degenerate, that is,
$\mathop{\rm rk}E=M$, then we obtain precisely Proposition 2.1 in
\cite{Pukh09b}. Proof of the latter proposition works in our case
without modifications.\vspace{0.3cm}

%%%%%%%%%%%%%%%%%%%%%%%%%%%%%%%%%%%%%%%%%%%%%%%%%%%%%%%%%%%%%%%%%
%%%%%%%%%%%%%%%%%%%%%%%%%%%%%%%%%%%%%%%%%%%%%%%%%%%%%%%%%%%%%%%%%
%%%%%%%%%%%%%%%%%%%%%%%%%%%%%%%   subsection 1.4

{\bf 1.4. Step 2: the germ of a quadratic singularity.} In this
subsection we consider $o\in W$ as a germ of a quadratic
singularity
$$
q_2(z_*)+q_3(z_*)+\dots =0,
$$
where $(z_*)=(z_1,\dots,z_M)$ (that is, disregarding the embedding
$W\subset{\mathbb P}^M$), so that $\mathop{\rm dim}W=M-1$. Let
$\varphi\colon W^+\to W$ be the blow up of the point $o$ and
$E=\varphi^{-1}(o)\subset W^+$ the exceptional divisor, a quadric
of rank $\mathop{\rm rk}q_2$ in ${\mathbb P}^{M-1}$. Consider an
effective divisor $D\ni o$ and assume that for the pair
$(W,\frac{1}{m}D)$ the point $o$ is an isolated centre of a
non-canonical singularity. Let $D^+\subset W^+$ be the strict
transform of the divisor $D$, so that $D^+=\varphi^*D-lE$ for some
$l\geq1$. Assume that $l\leq 2m$, so that the pair
$(W^+,\frac{1}{m}D^+)$ is not log canonical. Finally, let
$S\subset E$ be the centre of a non log canonical singularity of
that pair, which has the maximal dimension, in particular, $S$ is
not strictly contained in the centre of another non log canonical
singularity, if they exist. Obviously, the inequality
\begin{equation}\label{30.09.1}
\mathop{\rm mult}\nolimits_SD^+>m
\end{equation}
holds.\vspace{0.1cm}

The following claim generalizes Proposition 2.2 in \cite{Pukh09b}:
\vspace{0.1cm}

{\bf Proposition 1.2.} {\it Assume that the inequality
$$
\mathop{\rm dim}S+\mathop{\rm rk}q_2\geq M+3
$$
holds. Then one of the two cases takes place:\vspace{0.1cm}

(1) either $S$ is a hyperplane section of the quadric $E$ (the
simple case),\vspace{0.1cm}

(2) or there exists a hyperplane section $Z\supset S$ of the
quadric $E$, satisfying the inequality}
\begin{equation}\label{30.09.2}
\mathop{\rm mult}\nolimits_ZD^+>\frac{2m-l}{3}.
\end{equation}

{\bf Proof} is obtained partially by repeating the proof of
Proposition 2.2 in \cite{Pukh09b} word for word, partially by
reduction to that proposition via restricting the divisor $D$ onto
a generic section of the singularity $o\in W$ by a linear subspace
of dimension $\mathop{\rm rk}q_2$.\vspace{0.1cm}

More precisely, arguing as in \cite{Pukh09b}, we obtain from the
inequality (\ref{30.09.1}), that if $S\subset E$ is a prime
divisor, then $S\sim H_E$ is a hyperplane section of the quadric
$E$, that is, the case (1) takes place. Therefore, we assume that
$\mathop{\rm codim}(S\subset E)\geq 2$. Now, arguing as in
\cite{Pukh09b} (replacing Proposition 2.1 in that paper by
Proposition 1.1), we obtain that there exists a hyperplane section
$Z\supset S$, which is uniquely determined by the log pair
$(W^+,\frac{1}{m}D^+)$, satisfying the description of the case (2)
of Proposition 1.1.\vspace{0.1cm}

Now let us restrict the divisor $D$ onto the section
$W_{\Lambda}=W\cap\Lambda$ of the variety $W\subset{\mathbb C}^M$
by a generic linear subspace $\Lambda$ of dimension $\mathop{\rm
rk}q_2$. In this way we obtain the pair
$(W_{\Lambda},\frac{1}{m}D_{\Lambda})$ satisfying the assumptions
of Proposition 2.2 in \cite{Pukh09b} (the germ $o\in W_{\Lambda}$
is a germ of a non-degenerate quadratic singularity), the
subvariety $S_{\Lambda} =S\cap W^+_{\Lambda }$ is the centre of a
non log canonical singularity of the pair
$(W^+_{\Lambda},\frac{1}{m}D^+_{\Lambda})$, which has the maximal
dimension, so that the hyperplane section $Z_{\Lambda}=Z\cap
W^+_{\Lambda}$ satisfies the inequality
$$
\mathop{\rm
mult}\nolimits_{Z_{\Lambda}}D^+_{\Lambda}>\frac{2m-l}{3},
$$
which by genericity of the linear subspace $\Lambda$ implies the
required inequality (\ref{30.09.2}). Proposition 1.2 is shown.
Q.E.D. \vspace{0.3cm}

%%%%%%%%%%%%%%%%%%%%%%%%%%%%%%%%%%%%%%%%%%%%%%%%%%%%%%%%%%%%%%%%
%%%%%%%%%%%%%%%%%%%%%%%%%%%%%%%%%%%%%%%%%%%%%%%%%%%%%%%%%%%%%%%%
%%%%%%%%%%%%%%%%%%%%%%%   subsection 1.5

{\bf 1.5. Step 3: exclusion of the maximal singularity.} Let us
come back to the hypersurface $W\subset{\mathbb P}^M$ of degree
$M$ with an isolated quadratic singularity $o\in W$ of rank $\leq
5$. Let $\varphi\colon W^+\to W$ be its blow up,
$E=\varphi^{-1}(o)$ the exceptional quadric. Consider an effective
divisor $D\sim mH$, where $H$ is the class of a hyperplane section
of $W$, and let $D^+\sim mH-\nu E$ be its strict transform on
$W^+$.\vspace{0.1cm}

{\bf Proposition 1.3.} {\it The inequality $\nu\leq\frac32m$
holds.}\vspace{0.1cm}

{\bf Proof.} Since the inequality to be shown is linear in $D$,
without loss of generality we assume that $D$ is a prime divisor.
Assume the converse: $\nu>\frac32 m$. To begin with, consider the
first hypertangent divisor $D_2=\{q_2|_W=0\}$. Since by the
regularity condition (R3) $q_3|_E\not\equiv 0$, we have $D^+_2\sim
2H-3E$, which implies that the divisor $D_2$ is
reduced.\vspace{0.1cm}

{\bf Lemma 1.3.} {\it The divisor $D_2$ is
irreducible.}\vspace{0.1cm}

{\bf Proof.} Assume the converse. Then $D_2=\Delta_1+\Delta_2$,
where $\Delta_{1,2}$ are distinct hyperplane sections. Since the
quadric $E$ is irreducible, $\Delta^+_i\sim H-\alpha_iE$, where
$\alpha_i\in\{0,1\}$, so that we have
$D^+_2=\Delta^+_1+\Delta^+_2\sim 2H-\alpha E$, where
$\alpha\in\{0,1,2\}$. The contradiction proves the lemma. Q.E.D.
\vspace{0.1cm}

Therefore, $D$ and $D_2$ are distinct prime divisors, so that the
scheme-theoretic intersection $Y=(D\circ D_2)$ is well defined and
satisfies the inequality
$$
\frac{\mathop{\rm mult}_o}{\mathop{\rm
deg}}Y\geq\frac32\frac{2\nu}{mM}>\frac{9}{2M}.
$$
Now, repeating the proof of Proposition 3.1 in \cite{Pukh09b} word
for word, we obtain a contradiction by means of the method of
hypertangent divisors. Proposition 1.3 is shown.
Q.E.D.\vspace{0.1cm}

Now let us complete the proof of Theorem 3. Assume that the point
$o$ is an isolated centre of a non-canonical singularity of the
pair $(W,\frac{1}{m}D)$. By linearity of the Noether-Fano
inequality we may assume that $D$ is a prime divisor. Since
$\nu\leq\frac32m$, the pair $(W^+,\frac{1}{m}D^+)$ is not log
canonical and a certain irreducible subvariety $S\subset E$ is the
centre of a non log canonical singularity of that pair. We assume
that $S$ has the maximal dimension among all centres of non log
canonical singularities of the pair
$(W^+,\frac{1}{m}D^+)$.\vspace{0.1cm}

{\bf Proposition 1.4.} {\it The subvariety $S$ has codimension at
least 2 in the exceptional quadric $E$.}\vspace{0.1cm}

{\bf Proof} repeats the proof of Proposition 3.2 in \cite{Pukh09b}
word for word. Following the scheme of arguments in Sec. 3.2 in
\cite{Pukh09b}, we conclude that the second case of Proposition
1.2 takes place: there is a hyperplane section $Z\supset S$ of the
exceptional quadric $E$, satisfying the inequality
(\ref{30.09.2}). Let $P\subset{\mathbb P}^M$ be the (unique)
hyperplane, cutting out $Z$ on $E$, that is, $W^+_P\cap E=Z$,
where $W_P=W\cap P$. Obviously, the prime divisors $W_P$ and $D$
are distinct, so that the effective cycle $D_P=(D\circ W_P)$ of
codimension 2 satisfies the inequality
\begin{equation}\label{30.09.3}
\mathop{\rm mult}\nolimits_oD_P\geq\mathop{\rm mult}\nolimits_oD+
2\mathop{\rm mult}\nolimits_ZD^+>\frac43(l+m)>\frac83m.
\end{equation}
Now consider the pair $(W_P,\frac{1}{m}D_P)$. Its strict transform
$(W^+_P,\frac{1}{m}D^+_P)$ is not log canonical. We may assume
that the inequality $\mathop{\rm mult}_oD_P\leq 4m$ holds,
otherwise we obtain a contradiction, repeating the proof of
Proposition 1.3 word for word. The subvariety $S$ is contained in
the maximal centre $S'$ of a non log canonical singularity of the
pair $(W^+_P,\frac{1}{m}D^+_P)$. It is easy to see that $S'\subset
E_P=Z$ (otherwise, $\mathop{\rm dim}\varphi(S')\geq 5$, so that,
as $\mathop{\rm dim}\mathop{\rm Sing}W_P\leq 1$, there is a curve
$C\subset\varphi(S')$, $C\cap\mathop{\rm Sing}W_P=\emptyset$,
satisfying the inequality $\mathop{\rm mult}_CD_P> m$, which is
impossible for $D_P\sim mH_P$). For simplicity of notations we
assume that $S'=S$.\vspace{0.1cm}

Applying Proposition 1.2 once again, we obtain that one of the
following two cases takes place:\vspace{0.1cm}

(1) either $S$ is a hyperplane section of the quadric
$E_P$,\vspace{0.1cm}

(2) or there is a hyperplane section $Z^*\supset S$ of the quadric
$E_P$, satisfying the inequality
\begin{equation}\label{30.09.4}
\mathop{\rm mult}\nolimits_{Z^*}D^+_P>\frac{2m-l^*}{3},
\end{equation}
where $D^+_P\sim mH_P-l^*E_P$. By the inequality (\ref{30.09.3}),
the integer $l^*$ satisfies the inequality $l^*>\frac43m$. Now,
repeating the arguments in the beginning of Sec. 3.3 in
\cite{Pukh09b} word for word (using the regularity condition (R3)
instead of the condition (R2.2) in \cite{Pukh09b}), we exclude the
case (1).\vspace{0.1cm}

Now let us consider the hardest case (2). Since we can not use the
strong regularity condition (R2.2) that was used in
\cite{Pukh09b}, we need to slightly modify the arguments of Sec.
3.3 in that paper; in particular, we have to assume that $M\geq
14$. Let $R\subset P={\mathbb P}^{M-1}$ be the unique hyperplane,
cutting out $Z^*$ on the exceptional quadric $E_P$, that is,
$W^+_R\cap E_P=Z^*$, where $W_R=W_P\cap R$. Since $W_R\sim
H_P-E_P$, the pair $(W^+_P,W_R)$ is canonical. Furthermore,
$\mathop{\rm mult}_oW_R=2<\frac83$, so that by linearity of the
inequality (\ref{30.09.3}) and linearity of the condition of non
log canonicity of the pair $(W^+_P,\frac{1}{m}D^+_P)$ at $S$, we
may assume that $D_P$ does not contain the hyperplane section
$W_R$ as a component (in other words, removing that component, we
only make the inequality (\ref{30.09.3}) and the log Noether-Fano
inequality stronger). Therefore, we can take the effective cycle
$D_R=(D_P\circ W_R)$ of codimension 2 on $W_P$, which satisfies
the inequality
\begin{equation}\label{30.09.5}
\mathop{\rm mult}\nolimits_oD_R\geq\mathop{\rm
mult}\nolimits_oD_P+2\mathop{\rm
mult}\nolimits_{Z^*}D^+_P>\frac{28}{9}m.
\end{equation}
Since by the regularity condition (R3) the quadric $q_2|_R=0$ is
irreducible and $q_3|_{R\cap\{q_2=0\}}\not\equiv 0$, we may repeat
the proof of Lemma 1.3 and conclude that the divisor $D_2|_R$ is
irreducible and has the multiplicity precisely 6 at the point $o$.
Therefore,
$$
\frac{\mathop{\rm mult}\nolimits_o}{\mathop{\rm
deg}}(D_2|_R)=\frac{3}{M}.
$$
Let $Y=Y_3$ be an irreducible component of the effective cycle
$D_R$ with the maximal value of the ratio $(\mathop{\rm
mult}_o/\mathop{\rm deg})$. We have
$$
\frac{\mathop{\rm mult}_o}{\mathop{\rm deg}}Y>\frac{28}{9M},
$$
so that $Y\neq D_2|_R$, that is, $Y\not\subset D_2$ and we may
take the effective cycle $(Y_3\circ D_2)$ of codimension 2 on
$W_R$ and 4 on $W$, respectively. At least one of its components
$Y_4$ satisfies the inequality
$$
\frac{\mathop{\rm mult}_o}{\mathop{\rm
deg}}Y_4>\frac32\cdot\frac{28}{9M}=\frac{14}{3M}.
$$
Now let us apply the technique of hypertangent divisors to the
variety $W_R$ at the point $o$, satisfying the regularity
condition. We obtain a sequence of irreducible subvarieties
$$
Y_4,Y_5,\dots,Y_{M-2},
$$
$\mathop{\rm dim}Y_i=M-1-i$, where the curve $Y_{M-2}$ satisfies
the inequality
$$
\frac{\mathop{\rm mult}_o}{\mathop{\rm
deg}}Y_{M-2}>\frac{14}{3M}\cdot\frac54\cdot\frac65\cdot\dots\cdot
\frac{M-2}{M-3}=\frac{7(M-2)}{6M}.
$$
This is impossible for $M\geq 14$.\vspace{0.1cm}

Proof of Theorem 3 is complete. Q.E.D.\vspace{0.3cm}

%%%%%%%%%%%%%%%%%%%%%%%%%%%%%%%%%%%%%%%%%%%%%%%%%%%%%%%%%%%%%%%%%
%%%%%%%%%%%%%%%%%%%%%%%%%%%%%%%%%%%%%%%%%%%%%%%%%%%%%%%%%%%%%%%%%
%%%%%%%%%%%%%%%%%%%%%%   subsection 1.6

{\bf 1.6. Regular Fano hypersurfaces.} Let us consider Theorem 2.
The conditions (R2) and (R3) are checked by a routine dimension
count, which we skip. The condition (R4) checks quite elementary
(see Sec. 7.2 in \cite{Pukh10}) and we leave this work to the
reader. Let us consider the condition (R1) and outline the
proof.\vspace{0.1cm}

Let ${\cal P}$ be the linear space, consisting of tuples of
homogeneous polynomials $(p_1,\dots,p_N)$ of degrees $\mathop{\rm
deg} p_i=i+1$ on the projective space ${\mathbb P}^N$. Consider
the closed subset
$$
{\cal P}_{\rm non-reg}=\{\,(p_*)\in {\cal P}\,|\, \dim \{p_1=\dots
=p_N=0\}\geq 1\}.
$$

{\bf Proposition 1.5.} {\it The following equality holds:}
$$
\mathop{\rm codim} ({\cal P}_{\rm non-reg}\subset {\cal
P})=\frac{N(N+1)}{2}+2.
$$

{\bf Remark 1.2.} The claim can be made more precise: the closed
subset ${\cal P}_{\rm non-reg}$ is reducible and only one of its
components has the codimension given above, namely, the component,
consisting of such tuples $(p_*)$ that the closed subset
$\{p_1=\dots=p_N=0\}$ contains a line in ${\mathbb P}^N$. The
codimensions of the other components of the set ${\cal P}_{\rm
non-reg}$ are higher. However, we do not need this more precise
claim.\vspace{0.1cm}

{\bf Proof} of Proposition 1.5 is obtained by means of the methods
of the papers \cite{Pukh98b,Pukh01} (see also \cite[Chapter
3]{Pukh13a}): it is completely similar to the arguments of
\cite[Section 1]{Pukh98b}, when regularity of the sequence of
polynomials $p_1,\dots,p_N$ is violated for the first time at the
$i$-th polynomial, $i=1,\dots,N-1$, and to the arguments of
\cite[Section 3]{Pukh01} in the case when regularity is for the
first time violated at the last step, that is, the set
$$
\{p_1=\dots=p_{N-1}=0\}\subset {\mathbb P}^N
$$
is a curve and the polynomial $p_N$ vanishes identically on one of
its components. The computations are quite elementary and are left
to the reader. Q.E.D.\vspace{0.1cm}

{\bf Corollary 1.1.} {\it A generic (in the Zariski sense)
hypersurface $V$ satisfies the condition (R1) for $M\geq
13$.}\vspace{0.1cm}

{\bf Proof.} In the notations of the condition (R1) it is
sufficient to consider the worst case $c=3$. Taking into account
the dimension of the Grassmanian of subspaces of codimension 3 in
${\mathbb C}^{M+1}$ and the fact that the point $o\in V$ is
arbitrary, by Proposition 1.5 we get that the hypersurface $V$
satisfies the condition (R1), if the inequality
$$
\frac{(M-4)(M-3)}{2}+2-3(M-2)-M\geq 1
$$
holds. It is easy to check that the latter inequality is true for
$M\geq 13$. Proof is complete.\vspace{0.1cm}

Proof of Theorem 2 is now complete.

%%%%%%%%%%%%%%%%%%%%%%%%%%%%%%%%%%%%%%%%%%%%%%%%%%%%%%%%%%%%%%%%%
%%%%%%%%%%%%%%%%%%%%%%%%%%%%%%%%%%%%%%%%%%%%%%%%%%%%%%%%%%%%%%%%%
%%%%%%%%%%%%%%%%%%%%%%%%%%%%%%%%%%%%%%%%%%%%%%%%%%%%%%%%%%%%%%%%%
%%%%%%%%%%%%%%%%%%%%%%%%%%%%%%%%%%%%%%%%%%%%%%%%%%%%%%%%%%%%%%%%%
%%%%%%%%%%%%%%%%%%%%%%%%%%%%%%   SECTION 2

\section{Subvarieties of codimension two}

In this section we prove Theorem 4: if $B$ is a maximal subvariety
of codimension two for the system $\Sigma$, then $B$ is a section
of the hypersurface $V$ by a linear subspace of codimension two.
The proof makes use of the cone technique, see \cite[Chapter
2]{Pukh13a}. The main idea of our arguments is to consider
two-dimensional cones, swept out by {\it secant lines} of the
subvariety $B$.\vspace{0.3cm}

{\bf 2.1. The secant space of the subvariety $B$.} Assume that the
inequality (\ref{24.09.1}) holds. We need to show that $B=P\cap
V$, where $P\subset{\mathbb P}$ is a linear subspace of
codimension two. If $B$ is contained in a hyperplane, $B\subset
\Pi$, then the claim of the theorem is almost obvious: the
hyperplane section $V_{\Pi}=V\cap\Pi$ is a factorial variety,
$\mathop{\rm Pic}V_{\Pi}={\mathbb Z}H_{\Pi}$, where $H_{\Pi}$ is
the class of a hyperplane section, so that $B\sim mH_{\Pi}$ on
$V_{\Pi}$ for some $m\geq 1$. The restriction $\Sigma_{\Pi}$ of
the linear system $\Sigma$ на $V_{\Pi}$ is a non-empty system of
divisors, $\Sigma_{\Pi}\subset|2nH_{\Pi}|$, whereas $\mathop{\rm
mult}_B\Sigma_{\Pi}>n$, that is, $B$ is a fixed component of the
system $\Sigma_{\Pi}$ of multiplicity ($\mathop{\rm
mult_B\Sigma_{\Pi}})$. This implies that $m=1$, so that
$B\in|H_{\Pi}|$ is a hyperplane section of the variety
$V_{\Pi}\subset\Pi$, which is what we need.\vspace{0.1cm}

Starting from this moment, we assume that $B$ is not contained in
a hyperplane, that is, $\langle B\rangle={\mathbb P}$. Let us show
that the inequality (\ref{24.09.1}) is impossible for the mobile
linear system $\Sigma\subset|2nH|$. In order to do this, we assume
that this inequality is true and show that this assumption leads
to a contradiction.\vspace{0.1cm}

Define the {\it secant space}
$$
\mathop{\bf Sec}(B)\subset B\times B\times{\mathbb P}
$$
as the closure of the set
$$
\mathop{\bf Sec}\nolimits^*(B)\subset (B\times
B\backslash\Delta_B)\times{\mathbb P}, \quad \mathop{\bf
Sec}\nolimits^*(B)=\{(p,q,r)\,|\,r\in[p,q]\},
$$
where $\Delta_B\subset B\times B$ is the diagonal, $[p,q]$ is the
line, connecting the {\it distinct} points $p,q$. Let $\pi_B$ and
$\pi_{\mathbb P}$ be the projections of the irreducible variety
$\mathop{\rm Sec}(B)$ onto $B\times B$ and ${\mathbb P}$,
respectively.\vspace{0.1cm}

{\bf Proposition 2.1.} {\it The projection $\pi_{\mathbb P}$ is
surjective.}\vspace{0.1cm}

{\bf Proof} is given below in Subsection 2.3.\vspace{0.1cm}

Proposition 2.1 implies that the image of the restriction of the
projection $\pi_{\mathbb P}$ onto the set $\mathop{\bf Sec}^*(B)$
contains an open subset in ${\mathbb P}$. In the sequel, speaking
about a point $x$ of general position in ${\mathbb P}$, we will
always mean, in particular, that $x\not\in V$, so that the
restriction of the projection from the point $x$ onto $V$ is a
finite morphism $V\to{\mathbb P}^M$. Let $\mathop{\bf Sec}(B,x)$
and $\mathop{\bf Sec}^*(B,x)$ be the fibres of the projection
$\pi_{\mathbb P}$ and its restriction onto $\mathop{\bf
Sec}^*(B,x)$ over a point of general position $x\in{\mathbb
P}$.\vspace{0.1cm}

Obviously, $\mathop{\bf Sec}(B,x)$ can be considered as a closed
subset in $B\times B$, invariant with respect to the involution
$\tau\colon(p,q)\mapsto(q,p)$, and $\mathop{\bf Sec}^*(B,x)$ as a
closed subset in $B\times B\backslash\Delta_B$, where for a
sufficiently general point $x\in{\mathbb P}$)
$$
\mathop{\bf Sec}(B,x)=\overline{\mathop{\bf Sec}\nolimits^*(B,x)}.
$$
We have $\mathop{\rm dim}\mathop{\bf Sec}(B)=2M-3$, so that
$\mathop{\rm dim}\mathop{\bf Sec}(B,x)=M-4$. Taking the sections
of the closed set $\mathop{\bf Sec}(B,x)$ by generic very ample
divisors on $B\times B$, we obtain for every irreducible component
of the set $\mathop{\bf Sec}(B,x)$ a dense family of curves
$\Gamma$, where a general curve $\Gamma$ of the family is not
$\tau$-invariant and does not meet any fixed closed subset of
codimension $\geq 2$ in $\mathop{\bf Sec}(B)$. Set
$$
C_+=\pi_1(\Gamma),\quad C_-=\pi_2(\Gamma),
$$
where $\pi_{1,2}\colon B\times B\to B$ are the projections onto
the first and second direct factors, respectively. By
construction, $C_-$ is contained in the cone with the vertex $x$
and the base $C^+$ and the other way round. The properties of that
cone (swept out by the lines $[p,q]\ni x$, for $p,q\in\Gamma$) and
of the curves $C_{\pm}$ can be made more precise.\vspace{0.1cm}

The following fact is true.\vspace{0.1cm}

{\bf Proposition 2.2.} {\it For some positive integers $d_C$ and
$d_R$ there is an algebraic family
$$
{\cal A}=\{A=(C,R,x)\}
$$
of triples $(C,R,x)$, where $C$ is an effective 1-cycle of degree
$d_C$ on ${\mathbb P}$, $R$ is an effective 1-cycle of degree
$d_R$ on ${\mathbb P}$ and $x\in{\mathbb P}$ is a point,
satisfying the following conditions:\vspace{0.1cm}

1) the projection $\pi_{\mathbb P}\colon{\cal A}\to{\mathbb P}$,
$$
\pi_{\mathbb P}\colon(C,R,x)\mapsto x,
$$
is dominant, that is, $\pi_{\mathbb P}({\cal A})$ contains a
non-empty Zariski open subset,\vspace{0.1cm}

2) $C=C_++C_-$, where $C_{\pm}$ are distinct irreducible curves,
$C_{\pm}\subset B$,\vspace{0.1cm}

3) for the cone $C(x)\subset{\mathbb P}$ with the vertex $x$ and
the base $C_+$ we have: $C_{\pm}$ are sections of the cone and the
equality
$$
(C(x)\circ V)=C_++C_-+R=C+R
$$
holds, where the effective 1-cycle $R$ does not contain $C_{\pm}$
as a component,\vspace{0.1cm}

4) for any point $p\in C_{\pm}$ we have $x\not\in
T_pC_{\pm}$,\vspace{0.1cm}

5) if $p\in C_{\pm}$ is a singular point of the curve $C_{\pm}$,
then $x\not\in T_pV$,\vspace{0.1cm}

6) for any point of intersection $p\in C_+\cap C_-$ the generator
$[p,x]$ of the cone $C(x)$ has at the point $p$ a simple tangency
with the hypersurface $V$:
$$
([p,x]\cdot V)_p=2,
$$
that is, $p\not\in R$,\vspace{0.1cm}

7) the components of the curves $R$ sweep out} $V$:
$$
\overline{\bigcup_{A\in{\cal A}}R}=V.
$$
\vspace{0.1cm}

{\bf Proof}, which makes use of the construction, immediately
preceding the statement of Proposition 2.2, is given below in
Subsection 2.2.\vspace{0.1cm}

Let us complete the {\bf proof of Theorem 4.} In the notations of
Proposition 2.2, consider an arbitrary irreducible component of
the residual curve $R$, which we for simplicity denote by the same
symbol. Let $D\in\Sigma$ be a generic divisor. By the property 7),
we may assume that $R\not\subset D$. Since $B\subset D$ and,
moreover, $\mathop{\rm mult}_BD>n$, the inequality
\begin{equation}\label{30.09.2013.6}
2n\mathop{\rm deg}R=(R\cdot D)\geq\sum_{p\in R\cap B}(R\cdot D)_p>
n\sum_{p\in R\wedge B}\mathop{\rm mult_p}R
\end{equation}
holds, where the last sum is taken over the usual {\it and
infinitely near} points of intersection of the curve $R$ and the
subvariety $B$: the set of those points is denoted by the symbol
$R\wedge B$ (see \cite{Pukh06b}).\vspace{0.1cm}

Since $C_{\pm}\subset B$ and, by the property 6), the curve $R$
does not contain the points of the intersection $C_+\cap C_-$, we
have
\begin{equation}\label{30.09.2013.7}
R\wedge B\supset (R\wedge C_+)\sqcup(R\wedge C_-),
\end{equation}
$\sqcup$ means a disjoint union. By the property 5), the curve $R$
meets $C_{\pm}$ at non-singular points of those curves. The
following lemma is a version of a very well known claim
\cite{Pukh06b}.\vspace{0.1cm}

{\bf Lemma 2.1.} {\it The following equality holds:}
$$
\sum_{p\in R\wedge C_{\pm}}\mathop{\rm
mult}\nolimits_pR=\mathop{\rm deg}R.
$$
\vspace{0.1cm}

In the last formula we mean {\it any} choice of the sign $+$ or
$-$.\vspace{0.1cm}

Now from the inequality (\ref{30.09.2013.6}), taking into account
(\ref{30.09.2013.7}), we obtain: $2n\mathop{\rm
deg}R>n(\mathop{\rm deg}R+\mathop{\rm deg}R)$, which is
impossible. Q.E.D. for Theorem 4.\vspace{0.1cm}

{\bf Remark 2.1.} Repeating the previous arguments word for word,
we exclude the possibility of {\it two} maximal subvarieties of
codimension two for the system $\Sigma$. Therefore, the section
$B=V\cap P$ is uniquely determined.\vspace{0.3cm}

%%%%%%%%%%%%%%%%%%%%%%%%%%%%%%%%%%%%%%%%%%%%%%%%%%%%%%%%%%%%%%%%%%
%%%%%%%%%%%%%%%%%%%%%%%%%%%%%%%%%%%%%%%%%%%%%%%%%%%%%%%%%%%%%%%%%%
%%%%%%%%%%%%%%%%%%%%%%%%%   subsection 2.2

{\bf 2.2. Proof of technical facts.} Let us show Lemma 2.1. By
genericity of the curve $C=C_++C_-$, each of the curves $C_{\pm}$
is a section of the cone $C(x)$. The normalizations
$\widetilde{C_{\pm}}$ of these curves are naturally isomorphic.
Let $C^+(x)$ be the blow up of the vertex of the cone $C(x)$ and
$$
\widetilde{C}(x)=\widetilde{C}_{\pm}\times_{C_{\pm}}C^+(x)
$$
the non-singular ruled surface over $\widetilde{C}_{\pm}$, where
the smooth curves $\widetilde{C}_+$ and $\widetilde{C}_-$ are
realized as its sections. Set $\widetilde{R}$ to be the strict
transform of $R$ on $\widetilde{C}(x)$. By the properties 4)-6) at
each point $p\in(R\cap C_+)\sqcup(R\cap C_-)$ the corresponding
curve $C_+$ or $C_-$ is non-singular and transversal to the
generator of the cone $[p,x]$, so that
$$
\sum_{q\in R\wedge C_{\pm}(p)}\mathop{\rm mult}\nolimits_q R
=\sum_{q\in\widetilde{R}\wedge\widetilde{C}_{\pm}(p)} \mathop{\rm
mult}\nolimits_q\widetilde{R}=(\widetilde{R}\cdot\widetilde{C}_{\pm})_p,
$$
where the subset $R\wedge C_{\pm}(p)\subset R\wedge C_{\pm}$
consists of the point $p$ and infinitely near points over it,
where the point of the surface $\widetilde{C}(x)$, corresponding
to the point $p$ of intersection of the curves $R$ and $C_+$ (or
$C_-$), is denoted by the same symbol $p$. Therefore,
$$
\sum_{p\in R\wedge C_{\pm}}\mathop{\rm
mult}\nolimits_pR=\sum_{p\in\widetilde{R}\cap C_{\pm}}
(\widetilde{R}\cdot\widetilde{C}_{\pm})_p=
(\widetilde{R}\cdot\widetilde{C}_{\pm}),
$$
but the last number is equal to $\mathop{\rm deg}R$, see
\cite{Pukh06b}. Q.E.D. for the lemma.\vspace{0.1cm}

{\bf Proof of Proposition 2.2.} The construction, immediately
preceding the statement of Proposition 2.2, gives an algebraic
family ${\cal A}$, satisfying the property 1) by Proposition 2.1.
Let us show that, somewhat shrinking the family ${\cal A}$ (that
is, taking a Zariski open subset in that family), one can ensure
that the remaining properties 2)-7) hold. Indeed, our construction
yields in the general case distinct irreducible curves
$C_{\pm}\neq C_-$, so that the property 2) can be assumed. Proof
of Proposition 2.1 implies easily that a generic secant line
$[p,q]$ of the variety $B$ is not a 3-secant, that is, $C_{\pm}$
are sections of the cone $C(x)$, whereas $C_{\pm}$ come into the
1-cycle $(C(x)\circ V)$ with multiplicity 1, which gives the
property 3).\vspace{0.1cm}

For any point $p\in B$ we have $B\not\subset T_pV$ (since $B$ is
not contained in a hyperplane by assumption), and for that reason
for a general point $x\in{\mathbb P}$ the direction of the line
$[p,x]$ defines a field of directions on a proper closed subset of
the set $\mathop{\bf Sec}(B,x)$ (consisting of the points $p\in
B$, at which $[p,x]\subset T_pB$) and for that reason for a
general curve $\Gamma$ its projections $C_{\pm}$ are nowhere
tangent to the lines $[p,x]$ (that is, at no point $p\in
C_{\pm}$), that is, the property 4) is satisfied (one should also
take into account that for a general point $p\in B$ the set
$$
\mathop{\rm Sing}B\cap T_pV
$$
has dimension at most $M-4$, since the hypersurface $V$ is
non-singular and the linear system
$\sum\limits^{M+1}_{i=0}\lambda_i(\partial F/\partial x_i)$
defines a finite morphism ${\mathbb P}\to{\mathbb
P}$).\vspace{0.1cm}

The property 5) again follows from the fact that $B$ is not
contained in a hyperplane: obviously,
$$
\pi_B(\mathop{\rm Sing}(\mathop{\bf Sec}(B)))\subset(\mathop{\rm
Sing} B\times B)\cup(B\times\mathop{\rm Sing}B),
$$
so that the point $(p,q)\in\Gamma$ is a singularity of that curve
(for a general curve $\Gamma$) if and only if $p$ or $q$ belongs
to $\mathop{\rm Sing}B$. Since the pairs of points $(p,q)\in
B\times B$ such that $p\in\mathop{\rm Sing}B$ and $q\in T_pV$,
form a subset of codimension al least 2, a general curve $\Gamma$
does not contain such pairs. This proves the property
5).\vspace{0.1cm}

Let us consider the property 6). The subset $\pi^{-1}_B(\Delta_B)$
is a closed subset of codimension 1 in $\mathop{\bf Sec}(B)$,
which may consist of several irreducible components of different
codimensions. A general curve $\Gamma$ does not intersect the
components of codimension two, so we are only interested in the
divisorial components.\vspace{0.1cm}

It is easy to see that the closure of the set
$$
\pi^{-1}_B(\Delta_B\backslash\mathop{\rm
Sing}\Delta_B)\subset\mathop{\bf Sec}(B)
$$
is a prime Weil divisor on $\mathop{\bf Sec}(B)$. For a
non-singular point $p\in B\backslash\mathop{\rm Sing}B$ we have
$$
\pi^{-1}_B((p,p))=\{(p,p)\}\times T_pB.
$$
Let $(z_1,\dots,z_{M+1})$ be a system of affine coordinates with
the origin at the point $p\in{\mathbb P}$ and
$$
f=q_1+q_2+\dots+q_M
$$
the equation of the hypersurface $V$. By the condition on the rank
of the quadratic form $q_2$ we have $q_2|_{T_pB}\not\equiv 0$, so
that the set of triples $(p,p,x)\in\mathop{\bf Sec}(B)$ such that
\begin{equation}\label{30.09.2013.8}
([p,x]\cdot V)_p\geq 3
\end{equation}
has in $\mathop{\bf Sec}(B)$ codimension 2, which is what we need.
Therefore, it is sufficient to prove the property 6) for singular
points $p\in\mathop{\rm Sing}B$, that is $(p,p)\in\mathop{\rm
Sing}\Delta_B$.\vspace{0.1cm}

Let $Y\subset\mathop{\rm Sing}B$ be an irreducible subset of
codimension $\geq 2$ with respect to $B$. Since obviously for
$p\in Y$
$$
\pi^{-1}_B((p,p))\subset\{(p,p)\}\times T_pV
$$
and $q_2|_{T_pV}\not\equiv 0$, we obtain once again, that the set
of triples $(p,p,x)\in\mathop{\bf Sec}(B)$, such that the
inequality (\ref{30.09.2013.8}) holds, is of codimension at least
two in $\mathop{\bf Sec}(B)$, which is what we need. Therefore, it
is sufficient to consider a divisorial component
$Q\subset\mathop{\rm Sing}(B)$, $\mathop{\rm codim}(Q\subset
B)=1$.\vspace{0.1cm}

{\bf Lemma 2.2.}{\it There exists a non-empty Zariski open subset
$U_Q\subset Q$ such that for any point $p\in U_Q$ the set
$$
\pi_{\mathbb P}(\pi^{-1}_B((p,p)))\subset{\mathbb P}
$$
is a union of finitely many linear subspaces of dimension $M-1$,
contained in $T_pV$ and containing} $T_pQ$.\vspace{0.1cm}

{\bf Proof:} straightforward local computations. Since $Q$ is a
{\it divisorial} component of the set of singular points
$\mathop{\rm Sing}B$, over a non-empty Zariski open subset
$U_B\subset B$ with a non-empty intersection $U_Q=U_B\cap Q$, the
resolution of singularities of the variety $B$ is just the
normalization $\widetilde{B}\to B$, so that at every point $p\in
U_Q$ the variety $B$ admits a simple analytic parametrization and
easy local computations, which we omit, give an explicit
description of the limit set of secant lines $[q,r]$ when $q\to p$
and $r\to p$. Proof of the lemma is complete.\vspace{0.1cm}

Now for any point $p\in U_Q$ and some linear subspace
$\Pi\subset\pi_{\mathbb P}(\pi^{-1}_B((p,p)))$ we have
$q_2|_{\Pi}\not\equiv 0$, so that the closed set
$$
\{q_2=0\}\cap\pi_{\mathbb P}(\pi^{-1}_B((p,p)))\subset {\mathbb P}
$$
is of dimension $M-2$. Therefore, the set of triples
$(p,p,x)\in\mathop{\bf Sec}(B)$, satisfying the inequality
(\ref{30.09.2013.8}), where $p\in U_Q$, has the dimension
$(M-3)+(M-2)=2M-5$, that is, the codimension 2 in $\mathop{\bf
Sec}(B)$. This completes the proof of the property 6).
\vspace{0.1cm}

Finally, the property 7) is obvious (for instance, follows
immediately from the proof of Proposition 2.1, given
below).\vspace{0.1cm}

Proposition 2.2 is shown. Q.E.D.\vspace{0.3cm}

%%%%%%%%%%%%%%%%%%%%%%%%%%%%%%%%%%%%%%%%%%%%%%%%%%%%%%%%%%%%%%%%%%%%%
%%%%%%%%%%%%%%%%%%%%%%%%%%%%%%%%%%%%%%%%%%%%%%%%%%%%%%%%%%%%%%%%%%%%%
%%%%%%%%%%%%%%%%%%%%%%%%%   subsection 2.3

{\bf 2.3. The secant variety.} Set
$$
\mathop{\rm Sec}(B)=\pi_{\mathbb P}(\mathop{\bf Sec} (B)
)=\overline{\bigcup_{B\ni p\neq q\in B}[p,q]}={\mathbb P}
$$
and let us call $\mathop{\rm Sec}(B)$ {\it the secant variety} of
the subvariety $B\subset{\mathbb P}$ (as opposed to the secant
space, introduced in Subsection 2.1). We need to show that
$\mathop{\rm Sec}(B)={\mathbb P}$. Let $\alpha\colon{\mathbb
C}^{M+2}\backslash\{0\}\to{\mathbb P}$ be the canonical
projection. For a closed set $Y\subset{\mathbb P}$ the symbol
$Y^{\rm aff}$ stands for the affine cone
$$
\overline{\alpha^{-1}(Y)}=\alpha^{-1}(Y)\cup\{0\}\subset{\mathbb
C}^{M+2}.
$$
Let $\sigma\colon B^{\rm aff}\times B^{\rm aff}\times{\mathbb
C}^2\to{\mathbb C}^{M+2}$ be the map of taking the linear
combination
$$
\sigma\colon(v,w,(\lambda,\mu))\mapsto\lambda v+\mu w.
$$
Obviously, $\mathop{\rm Sec}(B)^{\rm aff}$ is the closure of the
image of the map $\sigma$. Furthermore, it is obvious that for a
non-singular point $p\in B$ the tangent space $T_vB^{\rm aff}$
does not depend on the choice of a non-zero vector
$v\in\alpha^{-1}(p)$ and for that reason we denote it by the
symbol $T_pB^{\rm aff}$. It is clear that the {\it embedded}
tangent space $T_pB\subset{\mathbb P}$ satisfies the equality
$$
(T_pB)^{\rm aff}=T_pB^{\rm aff}.
$$
Let $p,q\in B$ be a pair of non-singular points. Obviously, the
differential $d\sigma$ at the point $(v_0,w_0,(\lambda_0,\mu_0))$
is
$$
d\sigma\colon T_pB^{\rm aff}\times T_qB^{\rm aff}\times{\mathbb
C}^2\to{\mathbb C}^{M+2},
$$
$$
d\sigma\colon(v,w,(\lambda,\mu))\mapsto\lambda_0 v+\mu_0 w+\lambda
v_0+\mu w_0,
$$
so that for a non-singular point $r\in\mathop{\rm Sec}(B)$,
$r\in[p,q]$, we have $T_r\mathop{\rm Sec}(B)^{\rm aff}=\mathop{\rm
Im}d\sigma=T_pB^{\rm aff}+T_qB^{\rm aff}$ (taking into account
that $v_0\in T_pB^{\rm aff}$ and $w_0\in T_qB^{\rm aff}$). Set
$$
T(p,q)=T_pB^{\rm aff}+T_qB^{\rm aff}.
$$
Therefore, $\mathop{\rm Sec}(B)={\mathbb P}$ if and only if
$d\sigma$ is surjective, that is $T(p,q)={\mathbb
C}^{M+2}$.\vspace{0.1cm}

Assume now that $\mathop{\rm Sec}(B)\neq{\mathbb P}$ is a proper
irreducible subvariety. Since $\mathop{\rm codim}(B\subset{\mathbb
P})=3$, this implies that
$$
\mathop{\rm codim}(T_pB^{\rm aff}\cap T_qB^{\rm aff})\leq 5,
$$
and the latter holds for any non-singular points $p,q\in B$. Let
us show that our assumption leads to a
contradiction.\vspace{0.1cm}

The symbol $\pi_p$ stands for the linear projection ${\mathbb
P}\dashrightarrow{\mathbb P}^2$ from the tangent space $T_pB$ for
a non-singular point $p\in B$. The projection $\pi_p$ is the
projectivization of the linear map
$$
\pi_p^{\rm aff}\colon{\mathbb C}^{M+2}\to ({\mathbb
C}^{M+2}/T_pB^{\rm aff})\cong{\mathbb C}^3.
$$
The differential of the restriction of the latter map onto $B^{\rm
aff}$ is not surjective at a point of general position. Indeed,
for any smooth point $q\in B$ we have:
$$
\mathop{\rm dim}\pi^{\rm aff}_p(T_qB^{\rm aff})\leq 2.
$$
Therefore, $\pi_p(B)\neq{\mathbb P}^2$ and for that reason
$\pi_p(B)$ is either a point or some irreducible curve
$C\subset{\mathbb P}^2$. If $\pi_p(B)$ is a point or $C$ is a
line, then the subvariety $B$ is contained in a hyperplane, which
contradicts the assumption. Therefore, $\pi_p(B)=C$ is a curve of
degree $d\geq 2$.

Let $c\in C$ be a point of general position,
$$
B_c=\overline{(B\cap\pi^{-1}_p(c))\backslash T_pB}
$$
the fibre of the projection $\pi_p|_B$. Obviously, $B_c$ is a
closed subset of pure codimension two in the fibre
$\overline{\pi^{-1}_p(c)}\cong{\mathbb P}^{M-1}$. For that reason
the secant variety $\mathop{\rm Sec}(B_c)$ coincides with its
linear span $\langle B_c\rangle$ (it is sufficient to check this
almost obvious fact for a curve in ${\mathbb P}^3$).\vspace{0.1cm}

Therefore, we have three options:\vspace{0.1cm}

(1) $\mathop{\rm Sec}(B_c)={\mathbb P}^{M-1}$,\vspace{0.1cm}

(2) $\mathop{\rm Sec}(B_c)$ is a hyperplane in
$\overline{\pi^{-1}_p(c)}\cong{\mathbb P}^{M-1}$,\vspace{0.1cm}

(3) $\mathop{\rm Sec}(B_c)=B_c$ is a subspace of codimension two
in ${\mathbb P}^{M-1}$.\vspace{0.1cm}

Assume that the case (1) takes place. Since $\mathop{\rm
Sec}(B_c)\subset\mathop{\rm Sec}(B)$, we have
$$
\overline{\pi^{-1}_p(C)}\subset\mathop{\rm Sec}(B) .
$$
On the left we have an irreducible divisor in ${\mathbb P}$, so
that by our assumption that $\mathop{\rm Sec}(B)\neq{\mathbb P}$
the equality $\mathop{\rm Sec}(B)=\overline{\pi^{-1}_p(C)}$ holds.
However, it is obvious, that $\mathop{\rm Sec}(B)$ contains points
outside the set $\overline{\pi^{-1}_p(C)}$: let $c_1,c_2\in C$ be
a general pair of points, $q_i\in B_{c_i}$ general points, then
$[q_1,q_2]\subset\mathop{\rm Sec}(B)$, but
$\pi_p([q_1,q_2])=[c_1,c_2]\not\subset C$. This contradiction
excludes the case (1).\vspace{0.1cm}

The case (3) is impossible, as $V$ does not contain linear
subspaces of dimension $M-3$.\vspace{0.1cm}

Therefore, the case (2) takes place. Again we take a general pair
of points $c_1,c_2\in C$. Let $L=[c_1,c_2]\subset{\mathbb P}^2$ be
the line through them, $H=\overline{\pi^{-1}_p(L)}$ the
corresponding hyperplane in ${\mathbb P}$. Set also
$$
P_i=\overline{\pi^{-1}_p(c_i)}\subset H\cong{\mathbb P}^M.
$$
The linear space $P=T_pB$ is of codimension two in $H$ and
$P_1\cap P_2=P$. Furthermore, set
$$
B_i=B_{c_1}\,\,\mbox{и}\,\,\Pi_i=\mathop{\rm Sec}(B_i)= \langle
B_i\rangle,
$$
these are hyperplanes in $P_i$, $i=1,2$.\vspace{0.1cm}

{\bf Proposition 2.3.} {\it The following equality holds:}
$$
\mathop{\rm Sec}(B_1\cup B_2)=H.
$$

It is clear that since the points $c_1,c_2$ are general,
Proposition 2.3 implies the equality $\mathop{\rm Sec}(B)={\mathbb
P}$, which contradicts the initial assumption and proves
Proposition 2.1.\vspace{0.1cm}

{\bf Proof of Proposition 2.3.} None of the irreducible components
of the sets $B_1,B_2$ is a cone. Let $\Lambda\subset H$ be a
5-dimensional subspace of general position, $Q_i=P_i\cap\Lambda$
and $S_i=B_i\cap\Lambda$, $i=1,2$. Now $S_1,S_2$ are (possibly
reducible) surfaces in $\Lambda\cong{\mathbb P}^5$, the linear
spans $\langle S_i\rangle$ of which are 3-planes
$R_i=\Pi_i\cap\Lambda$. The components of the surfaces $S_1,S_2$
are not cones and for that reason
$$
\bigcap_{s\in S_i}T_sS_i=\emptyset
$$
for $i=1,2$. If $R_{12}=R_1\cap R_2$ is a line, then we conclude
that for a general pair of points $(s_1,s_2)\in S_1\times S_2$ the
planes $T_{s_1}S_1$ and $T_{s_2}S_2$ are disjoint. This implies,
that $\mathop{\rm Sec}(S_1\cup S_2)=\Lambda$, so that Proposition
2.3 is shown in this case.\vspace{0.1cm}

Therefore we assume that $R_{12}$ is a plane, that is,
$$
R_1\cap R_2=R_1\cap Q_2=Q_1\cap R_2.
$$
By the genericity of the subspace $\Lambda$ this means that
$$
\Pi_{12}=\Pi_1\cap\Pi_2=\Pi_1\cap P_2=P_1\cap\Pi_2,
$$
and for that reason $\Pi_{12}=\Pi_1\cap P=\Pi_2\cap P$. The points
$c_1,c_2$ are chosen independently of each other, so that we can
conclude that there exists (a uniquely determined) hyperplane
$Q\subset P$ such that for a point of general position $c\in C$ we
have
$$
\langle B_c\rangle =\mathop{\rm Sec}(B_c)\supset Q.
$$
Let $\pi_Q\colon{\mathbb P}\dashrightarrow{\mathbb P}^3$ be the
projection from the linear subspace $Q$. By what we proved,
$\pi_Q(\langle B_c\rangle)$ is a point and for that reason
$\pi_Q(B_c)$ is a point, so that the image $\pi_Q(B)$ is a curve
$C^+$ (the projection of which from the point $\pi_Q(P)$ is the
curve $C\subset{\mathbb P}^2$). If $C^+$ is contained in some
plane in ${\mathbb P}^3$, then $B$ is contained in some hyperplane
in ${\mathbb P}$, either, which contradicts our assumption. Thus
$$
\langle C^+\rangle =\mathop{\rm Sec}(C^+)={\mathbb P}^3.
$$
Now let $\xi_1,\xi_2\in C^+$ be a general pair of points,
$\Lambda_i=\pi^{-1}_Q(\xi_i)\subset{\mathbb P}$ the corresponding
subspaces of codimension 3, $B^+_i=\overline{\pi^{-1}_Q(\xi_i)\cap
B}$ the fibres of the projection $\pi_Q|_B$. We know that
$B^+_i\subset\Lambda_i$ are hypersurfaces (possibly reducible) and
$$
\langle B^+_i \rangle =\mathop{\rm Sec}(B^+_i)=\Lambda_i.
$$
Since $B^+_i$ are not cones, we conclude that
$$
\mathop{\rm Sec}(B^+_i\cup B^+_2)=\pi^{-1}_Q([\xi_1,\xi_2]),
$$
whence, finally, it follows that  $\mathop{\rm Sec}(B)={\mathbb
P}$. Proof of Propositions 2.3 and 2.1 is complete.

%%%%%%%%%%%%%%%%%%%%%%%%%%%%%%%%%%%%%%%%%%%%%%%%%%%%%%%%%%%%%%%%%%%
%%%%%%%%%%%%%%%%%%%%%%%%%%%%%%%%%%%%%%%%%%%%%%%%%%%%%%%%%%%%%%%%%%%
%%%%%%%%%%%%%%%%%%%%%%%%%%%%%%%%%%%%%%%%%%%%%%%%%%%%%%%%%%%%%%%%%%%
%%%%%%%%%%%%%%%%%%%%%%%%%%%%%%%%%%%%%%%%%%%%%%%%%%%%%%%%%%%%%%%%%%%
%%%%%%%%%%%%%%%%%%%%%   SECTION 3

\section{Infinitely near case. I. Preparatory work}

In this section we start the proof of Theorem 5, that is, the
exclusion of the infinitely near case. Here we carry out
preparatory work: we come over to a hyperplane section of the
hypersurface $V$, in order to use the $8n^2$-неравенство, list all
particular cases that need to be considered and obtain aprioric
estimates for the multiplicity of the self-intersection. We use
the following tools: the inversion of adjunction, the technique of
counting multiplicities and the method of hypertangent
divisors.\vspace{0.3cm}

{\bf 3.1. The method of hypertangent divisors.} Let $\Sigma\subset
|2nH|$ be a mobile linear system with no maximal subvarieties of
codimension two. Fix a maximal singularity
$E^*\subset\widetilde{V}$ of the system $\Sigma$ with the centre
$B\subset V$ of maximal dimension.\vspace{0.1cm}

{\bf Lemma 3.1.} $B$ {\it is a point or a curve on
$V$.}\vspace{0.1cm}

{\bf Proof.} By the $4n^2$-inequality we have
$$
\mathop{\rm mult}\nolimits_BZ>4n^2,
$$
where $Z=(D_1\circ D_2)$ is the self-intersection of the system
$\Sigma$. Since $Z\sim 4n^2H^2$, by \cite[Proposition 5]{Pukh02b}
it follows that $\mathop{\rm dim}B\leq 1$. Q.E.D. for the
lemma.\vspace{0.1cm}

The cases $\mathop{\rm dim}B=1$ and $\mathop{\rm dim}B=0$ are
dealt with in word for word the same way, the assumption on the
existence of a maximal singularity leads to a contradiction,
excluding both cases. We will assume that $B=o$ is a point. The
following fact is true.\vspace{0.1cm}

{\bf Proposition 3.1.} {\it The following inequality holds:}
$\mathop{\rm mult}_o\Sigma\leq 3n$.\vspace{0.1cm}

{\bf Proof.} Assume the converse: $\mathop{\rm mult}_o\Sigma
>3n$. Let $T=T_oV\cap V$ be the intersection of the hypersurfaces $V$
with the tangent hyperplane. Obviously,
$$
T\subset T_oV\cong{\mathbb P}^M
$$
is a Fano hypersurface with the isolated double point $o\in T$.
The tangent cone at the point $o$ is the quadric
$\{q_2|_{\{q_1=0\}}=0\}$. For a generic divisor $D\in\Sigma$ we
have $D\neq T$, so that $(D\circ T)$ is an effective cycle of
codimension two, satisfying the inequality
$$
\mathop{\rm mult}\nolimits_o(D\circ T)>6n.
$$
Let $Y$ be a component of the cycle $(D\circ T)$ with the maximal
value of the ratio $(\mathop{\rm mult}_o/\mathop{\rm deg})$.
Therefore, the prime divisor $Y\subset T$ satisfies the inequality
$$
\frac{\mathop{\rm mult}_o}{\mathop{\rm deg}}Y>\frac{3}{M}.
$$
The first hypertangent divisor $T_2=\{q_2|_T=0\}$ is irreducible
and by the regularity conditions $\mathop{\rm mult}_oT_2=6$,
$\mathop{\rm deg}T_2=2M$, so that $T_2\neq Y$. Let us form the
effective cycle $(\{q_2|_T=0\}\circ Y)$ and choose in it an
irreducible component $Y_3$ with the maximal value of the ratio
$(\mathop{\rm mult}_o/\mathop{\rm deg})$. Now we apply to $Y_3$
the standard technique of hypertangent linear systems
\cite[Chapter 3]{Pukh13a}: take generic hypertangent divisors
$$
T_i\in|iH-(i+1)E|,\,\,i=4,\dots,M-1,
$$
and construct a sequence of irreducible subvarieties
$$
Y_3,Y_4,\dots,Y_{M-1},
$$
where $\mathop{\rm codim}(Y_i\subset V)=i$, and $Y_{i+1}$ is an
irreducible component of the effective cycle $(Y_i\circ T_{i+1})$
with the maximal value of $(\mathop{\rm mult}_o/\mathop{\rm
deg})$. For the curve $Y_{M-1}$ we have the inequality
\begin{equation}\label{11.05.2}
1\geq\frac{\mathop{\rm mult}_o}{\mathop{\rm deg}}Y_{M-1}>
\frac{3}{M}\cdot\frac32\cdot\frac54\cdot
\frac65\cdot\dots\cdot\frac{M}{M-1}=\frac98\cdot
\end{equation}
This contradiction proves our proposition. Q.E.D.\vspace{0.1cm}

Arguing in a similar way, we obtain the following
fact.\vspace{0.1cm}

{\bf Proposition 3.2.} {\it For any irreducible subvariety
$Y\subset V$ of codimension two the following inequality holds:}
\begin{equation}\label{1.10.2013.1}
\frac{\mathop{\rm mult}_o}{\mathop{\rm deg}}Y\leq\frac{3}{M}.
\end{equation}

{\bf Proof.} Set again $T=T_oV\cap V$. If $Y\subset T$, that is,
$Y$ is a prime divisor on the hypersurface $T\subset
T_oV\cong{\mathbb P}^M$, then we argue in word for word the same
way as in the proof of Proposition 3.1, deriving a contradiction
from the assumption that the inequality (\ref{1.10.2013.1}) does
not hold.\vspace{0.1cm}

Assume now, that $Y\not\subset T$ and the inequality
(\ref{1.10.2013.1}) is not true. Then for the effective cycle
$(Y\circ T)$ the inequality
$$
\frac{\mathop{\rm mult}_o}{\mathop{\rm deg}}(Y\circ T)>\frac{6}{M}
$$
holds, so that there is an irreducible component $Y_3$ of the
cycle $(Y\circ T)$, satisfying that inequality. Now we argue in
the same way as in the proof of Proposition 3.1: taking generic
hypertangent divisors $T_4,\dots,T_{M-1}$, we construct a sequence
of irreducible subvarieties $Y_3,Y_4,\dots,Y_{M-1}$, the last one
of which is a curve, satisfying the inequality
$$
1\geq\frac{\mathop{\rm mult}_o}{\mathop{\rm
deg}}Y_{M-1}>\frac{6}{M}
\cdot\frac54\cdot\frac65\cdot\dots\cdot\frac{M}{M-1}=\frac64.
$$
This contradiction completes the proof of the
proposition.\vspace{0.1cm}

{\bf Remark 3.1.} The second part of the proof of Proposition 3.2
gives a much stronger estimate for the ratio $(\mathop{\rm
mult}_o/\mathop{\rm deg})$ in the case $Y\not\subset T$:
$$
\frac{\mathop{\rm mult}_o}{\mathop{\rm deg}}Y\leq\frac{2}{M}.
$$
\vspace{0.1cm}

Now let $\Pi\subset{\mathbb P}$ be a linear subspace of
codimension 1 or 2, containing the point $o$, but not contained in
the hyperplane $T_oV$, so that $V_{\Pi}=V\cap\Pi$ is an
irreducible hypersurface of degree $M$ in $\Pi$, non-singular at
the point $o$. Let $Y\subset V_{\Pi}$ be an irreducible subvariety
of codimension two.\vspace{0.1cm}

{\bf Proposition 3.3.} (i) {\it If $Y\subset T_oV_{\Pi}$, then for
$\mathop{\rm codim} (\Pi\subset{\mathbb P})=j\in \{1,2\}$ the
following estimate holds:}
$$
\frac{\mathop{\rm mult}_o}{\mathop{\rm deg}}Y\leq\frac{1}{M}
\max\left( 3,\frac{8M}{3(M-j)}\right).
$$

(ii) {\it If $Y\not\subset T_oV_{\Pi}$, then for $\mathop{\rm
codim} (\Pi\subset{\mathbb P})=j\in \{1,2\}$ the following
estimate holds}
$$
\frac{\mathop{\rm mult}_o}{\mathop{\rm deg}}Y\leq\frac{1}{M}
\left(\frac{2M}{M-j}\right).
$$
\vspace{0.1cm}

(The somewhat strange writing of the right hand part of the
inequality (ii) will become clear below.)\vspace{0.1cm}

{\bf Proof.} (i) Repeating the arguments of the first part of the
proof of Proposition 3.2 word for word and taking into account the
regularity conditions for the hypersurface $V_{\Pi}$, we obtain
the inequality
$$
1>\frac{3}{M}\cdot\frac32\cdot\frac54\cdot\frac65\cdot\dots
\cdot\frac{M-1}{M-2}=\frac98\left(1-\frac{1}{M}\right)
$$
for $\mathop{\rm codim}(\Pi\subset{\mathbb P})=1$, and the
inequality
$$
1>\frac{3}{M}\cdot\frac32\cdot\frac54\cdot\frac65\cdot\dots
\cdot\frac{M-2}{M-3}=\frac98\left(1-\frac{2}{M}\right)
$$
for $\mathop{\rm codim}(\Pi\subset{\mathbb P})=2$. It is easy to
see that these inequalities are impossible. The contradiction
proves the claim (i).\vspace{0.1cm}

In the case (ii) we repeat the arguments of the second part of the
proof of Proposition 3.2 word for word, once again taking into
account that by the considerations of dimension we take
$\mathop{\rm codim}(\Pi\subset{\mathbb P})=1$ or 2 hypertangent
divisors less. Again we get a contradiction, which proves the
claim (ii).\vspace{0.1cm}

Proof of Proposition 3.3 is complete.\vspace{0.3cm}

%%%%%%%%%%%%%%%%%%%%%%%%%%%%%%%%%%%%%%%%%%%%%%%%%%%%%%%%%%%%%%%%%%%%
%%%%%%%%%%%%%%%%%%%%%%%%%%%%%%%%%%%%%%%%%%%%%%%%%%%%%%%%%%%%%%%%%%%%
%%%%%%%%%%%%%%%%%%%%%%%   subsection 3.2

{\bf 3.2. The restriction onto a hyperplane section.} The next
step in the proof of Theorem 5 is the restriction of the linear
system $\Sigma$ onto a suitable hyperplane section of the variety
$V$, which allows us to make the estimate for the multiplicity of
the self-intersection at the point $o$ twice stronger. If the
inequality $\mathop{\rm mult}_oZ>8n^2$ holds, then this step can
be skipped, considering below instead of the hyperplane section
$P\ni o$ the hypersurface $V$ itself: in that case, the dimension
does not drop and all estimates become only stronger, so that the
proof given below works without any modifications. Keeping this in
mind, assume that $\mathop{\rm mult}_oZ\leq 8n^2$. The following
fact is true.\vspace{0.1cm}

{\bf Proposition 3.4 (the $8n^2$-inequality).} {\it There exists a
subspace $\Pi\subset E$ of codimension 2 (uniquely determined by
the system $\Sigma$), satisfying the inequality}
$$
\mathop{\rm mult}\nolimits_oZ+ \mathop{\rm
mult}\nolimits_{\Pi}Z^+>8n^2.
$$

{\bf Proof:} this is \cite[Proposition 4.1]{Pukh10}.\vspace{0.1cm}

Now let us consider the linear system $|H-\Pi|$, consisting of
hyperplane sections that cut out $\Pi$ on $E$, that is, for a
general divisor $P\in|H-\Pi|$ we have: $P\in |H|$ is a hyperplane
section, smooth at the point $o$ and $P^+\supset\Pi$. Obviously,
$$
\mathop{\rm dim}|H-\Pi|=2\quad \mbox{и}\quad \mathop{\rm codim}
\mathop{\rm Bs}|H-\Pi|=3.
$$
Therefore for a general divisor $P\in |H-\Pi|$ the effective cycle
$Z_P=(Z\circ P)$ of codimension two is well defined and satisfies
the inequality
$$
\mathop{\rm mult}\nolimits_oZ_P=\mathop{\rm mult}\nolimits_oZ+
\mathop{\rm mult}\nolimits_{\Pi}Z^+>8n^2.
$$
Let $\Sigma_P=\Sigma|_P$ be the restriction of the linear system
$\Sigma$ onto $P$. Obviously, $\Sigma_P\subset |2nH_P|$, where
$H_P=H|_P$ is the positive generator of the group $\mathop{\rm
Pic}P\cong{\mathbb Z}$, whereas the system $\Sigma_P$ is mobile
(has no fixed components). The cycle $Z_P$ is the
self-intersection of the system $\Sigma_P$:
$$
Z_P=(D_1\circ D_2),
$$
where $D_1,D_2\in\Sigma_P$ are generic divisors. The variety $P$
is a hypersurface of degree $M$ in ${\mathbb P}^M$, which may have
isolated singular points, but the point $o\in P$ itself is
non-singular.\vspace{0.1cm}

{\bf Proposition 3.5.} {\it The pair $(P,\frac{1}{n}\Sigma_P)$ is
not log canonical at the point $o$, that is, it has a non log
canonical singularity with the centre at that point. If the pair
$(P,\frac{1}{n}\Sigma_P)$ has a non-canonical singularity with the
centre $B\ni o$, $B\neq o$, then either $\mathop{\rm dim}B\leq 2$,
or $B=\Delta=\mathop{\rm Bs}|H-\Pi|$ (and in the latter case the
inequality $\mathop{\rm mult}\nolimits_{\Delta}\Sigma>n)$ holds}
.\vspace{0.1cm}

{\bf Proof.} The first claim follows from the inversion of
adjunction \cite{Kol93}. Let us consider the second one (it is not
used in the subsequent proof). If $B\neq\Delta$, then $\mathop{\rm
codim}(B\subset P)\geq 3$ (otherwise, by the genericity of $P$,
the original system $\Sigma$ has a maximal subvariety of
codimension 2, which is not true by assumption). Therefore, the
$4n^2$-inequality holds:
$$
\mathop{\rm mult}\nolimits_BZ_P>4n^2.
$$
Let $Q\in |H_P|$ be a general (in particular, everywhere
non-singular) hyperplane section of $P$ and $Z_Q=(Z_P\circ Q)$.
Then on $Q$ the cycle $Z_Q\sim 4n^2H_Q$ of codimension two
satisfies the inequality
$$
\mathop{\rm mult}\nolimits_{B\cap Q}Z_Q>4n^2
$$
and $\mathop{\rm dim}B\cap Q\leq 1$ by Proposition 5 in
\cite{Pukh02b}. Proof is complete.\vspace{0.1cm}

Note that by the genericity of the hyperplane section $P$ the
linear system $\Sigma_P$ satisfies the inequality
$$
\nu=\mathop{\rm mult}\nolimits_o\Sigma_P(=\mathop{\rm
mult}\nolimits_o\Sigma)\leq 3n.
$$
Now let $\Pi_1\subset\Pi_2$ be a generic pair of linear subspaces
of dimensions 5 and 5 in ${\mathbb P}^M=\langle P\rangle$,
containing the point $o$, and $X_i=P\cap\Pi_i$ the corresponding
sections of the hypersurface $P$. By the inversion of adjunction
the pair $(X_i,\frac{1}{n}\Sigma_i)$, where
$\Sigma_i=\Sigma_P|_{X_i}$, has the point $o$ as an isolated
centre of a non log canonical singularity. Let $X^+_i\subset P^+$
be the strict transform of $X_i$, so that
$$
\varphi_i\colon X^+_i\to X_i
$$
is the blow up of the point $o\in X_i$ and $E^{(i)}=E\cap X^+_i$
the exceptional divisor of the morphism $\varphi_i$. The pairs
\begin{equation}\label{22.05.1}
\Box_1=\left(X^+_1,\frac{1}{n}\Sigma^+_1+\frac{\nu-3n}{n}E^{(1)}\right)
\end{equation}
and
\begin{equation}\label{22.05.2}
\Box_2=\left(X^+_2,\frac{1}{n}\Sigma^+_2+\frac{\nu-4n}{n}E^{(2)}\right)
\end{equation}
are not log canonical (recall that $\nu\leq 3n$) and satisfy the
conditions of the connectedness principle with respect to the
birational morphisms $\varphi_1$ and $\varphi_2$, respectively
(see \cite[Section 17.4]{Kol93}). The centre of any non log
canonical singularity of the pair $\Box_i$, intersecting
$E^{(i)}$, is contained in $E^{(i)}$ (Proposition 3.5), so that we
conclude that the union $LCS(\Box)_i$ of centres of non log
canonical singularities of the pair $\Box_i$, intersecting
$E^{(i)}$, is a connected closed subset in $E^{(i)}$. Recall that
$E^{(1)}\cong{\mathbb P}^3$ and $E^{(2)}\cong{\mathbb P}^4$. For
the pair $\Box_1$ there are three options:\vspace{0.1cm}

--- $LCS(\Box_1)$ is a point $p\in E^{(1)}$,\vspace{0.1cm}

--- $LCS(\Box_1)$ is a connected curve,\vspace{0.1cm}

--- $LCS(\Box_1)$ is a union of curves and surfaces, and
in this union there is is at least one surface.\vspace{0.1cm}

For the pair $\Box_2$ there are, respectively, four options,
$\mathop{\rm dim}LCS(\Box_2)\in\{0,1,2,3\}$, and if $LCS(\Box_2)$
is zero-dimensional, then this set consists of one
point.\vspace{0.1cm}

Now looking at the pair
$$
\Box_{12}=\left(X^+_2,\frac{1}{n}\Sigma^+_2+
\frac{\nu-3n}{n}E^{(2)}\right),
$$
we see that $LCS(\Box_{12})$ is either a line in
$E^{(2)}\cong{\mathbb P}^4$, or a connected union of surfaces
(every hyperplane section of which is connected), or a union of
surfaces and divisors in $E^{(2)}$. Since the pair $\Box_{12}$ is
obviously ``more effective'' than the pair $\Box_2$, we have the
inclusion
$$
LCS(\Box_2)\subset LCS(\Box_{12}),
$$
in particular, $(LCS(\Box_2)\cap X^+_1)\subset
LCS(\Box_1)$.\vspace{0.1cm}

Now let us come back to the hypersurface $P$ and its blow up
$\varphi_P\colon P^+\to P$ at the point $o$. From what was said,
it follows that the pairs
$$
\Box=\left(P^+,\frac{1}{n}\Sigma^+_P+\frac{\nu-3n}{n}E_P\right)
$$
and
$$
\Box^*=\left(P^+,\frac{1}{n}\Sigma^+_P+\frac{\nu-4n}{n}E_P\right)
$$
are not log canonical, and moreover, one of the following six
cases takes place.\vspace{0.1cm}

{\bf Case 1.1.} There are non log canonical singularities of the
pairs $\Box^*$ and $\Box$, the centres of which on $P^+$ are
linear subspaces $\Theta\subset\Lambda\subset E_P$ of codimension
4 and 5, respectively.\vspace{0.1cm}

{\bf Case 1.2.} There exists a non log canonical singularity of
the pair $\Box^*$, the centre of which on $P^+$ is a linear
subspace $\Lambda\subset E_P$ of codimension 3.\vspace{0.1cm}

{\bf Case 2.1.} There exist non log canonical singularities of the
pairs $\Box^*$ and $\Box$, the centres of which on $P^+$ are a
linear subspace $\Theta\subset E_P$ of codimension 4 and an
irreducible subvariety $B\subset E_P$ of codimension 2,
respectively, where $\Theta\subset B$.\vspace{0.1cm}

{\bf Case 2.2.} There are non log canonical singularities of the
pairs $\Box^*$ and $\Box$, the centres of which on $P^+$ are
irreducible subvarieties $B^*\subset B\subset E_P$ of codimension
3 and 2, respectively.\vspace{0.1cm}

{\bf Case 2.3.} There is a non log canonical singularity of the
pair $\Box^*$, the centre of which on $P^+$ is an irreducible
subvariety $B\subset E_P$ of codimension 2.\vspace{0.1cm}

{\bf Case 3.} There is a non log canonical singularity of the pair
$\Box$, the centre of which on $P^+$ is an irreducible subvariety
$B\subset E_P$ of codimension 1.\vspace{0.1cm}

The six cases listed above correspond to three possible values of
the integer $\mathop{\rm dim}LCS(\Box_1)$, taking into account the
type of the set $LCS(\Box_2)$.\vspace{0.1cm}

The last case is the simplest one.\vspace{0.1cm}

{\bf Proposition 3.6.} {\it The case 3 does not realize:}
$\mathop{\rm codim}(B\subset E_P)\geq 2$.\vspace{0.1cm}

{\bf Proof.} Assume the converse: $B\subset E_P$ is a prime
divisor. We argue as in the proof of Proposition 4.1 in
\cite{Pukh10} or in \cite{Ch06b}: for the self-intersection $Z_P$
of the system $\Sigma_P$, taking into account that the pair $\Box$
is not log canonical at $B$, we obtain the estimate
$$
\mathop{\rm mult}\nolimits_oZ_P>
\nu^2+4(4-\frac{\nu}{n})n^2=(\nu-2n)^2+12n^2\geq 12n^2.
$$
Therefore, there is an irreducible subvariety $Y\subset P$ of
codimension two, satisfying the inequality
$$
\frac{\mathop{\rm mult}_o}{\mathop{\rm deg}}Y>\frac{3}{M}.
$$
However, this contradicts Proposition 3.3. Proposition 3.6 is
shown. Q.E.D.\vspace{0.1cm}

{\bf Remark 3.2.} Once again, we emphasize that Proposition 3.3
implies the inequality
$$
\mathop{\rm mult}\nolimits_oZ_P\leq 12n^2,
$$
which we will use in the sequel without special
references.\vspace{0.3cm}

%%%%%%%%%%%%%%%%%%%%%%%%%%%%%%%%%%%%%%%%%%%%%%%%%%%%%%%%%%
%%%%%%%%%%%%%%%%%%%%%%%%%%%%%%%%%%%%%%%%%%%%%%%%%%%%%%%%%%
%%%%%%%%%%%%%%%%%%   subsection 3.3

{\bf 3.3. The techniques of counting multiplicities: the aprioric
estimates.} Following the standard procedure of the method of
maximal singularities, let us obtain now bounds from below for the
multiplicities of the cycle $Z_P$, improving the
$8n^2$-inequality. We call these estimates aprioric, because they
do not make use the additional geometric information available in
the cases 1.1-2.3. To exclude those cases, the aprioric estimates
are not sufficient and we will need some additional work, which
will be carried out in Sections 4,5.\vspace{0.1cm}

{\bf Proposition 3.7.} (i) {\it If the case 1.1 takes place, then
the following inequalities hold:
\begin{equation}\label{04.06.1}
\mathop{\rm mult}\nolimits_oZ_P+ \mathop{\rm
mult}\nolimits_{\Lambda}Z^+_P>12n^2
\end{equation}
and} $\mathop{\rm mult}_{\Theta}Z_P>4n^2$. {\it If the case 1.2
takes place, then the following estimate holds:}
\begin{equation}\label{04.06.2}
\mathop{\rm mult}\nolimits_{\Lambda}Z_P>4n^2.
\end{equation}

(ii) {\it If either of the cases 2.1 or 2.2 takes place, then the
following inequality holds
\begin{equation}\label{29.06.1}
\mathop{\rm mult}\nolimits_oZ_P+ \mathop{\rm mult}\nolimits_
BZ^+_P>12n^2,
\end{equation}
in addition in the case 2.1 the estimate $\mathop{\rm
mult}_{\Theta}Z^+_P>4n^2$ and in the case 2.2 the estimate
$\mathop{\rm mult}_{B^*}Z^+_P>4n^2$ hold.}\vspace{0.1cm}

(iii) {\it If the case 2.3 takes place, then the following
inequality holds:}
$$
\mathop{\rm mult}\nolimits_BZ^+_P>4n^2.
$$

{\bf Proof.} All the inequalities, listed above, belong to one of
the two types: the type (\ref{04.06.1}) for a non log canonical
singularity of the pair $\Box$ and the type (\ref{04.06.2}) for a
singularity of the pair $\Box^*$. The proofs for each of the two
types are completely identical, and for this reason we will show
only these two inequalities.\vspace{0.1cm}

Let us prove the inequality (\ref{04.06.1}). It is true under a
weaker assumption that the pair $\Box$ has a {\it non canonical}
singularity, the centre of which is the subspace $\Lambda$. This
is what we will assume. Let
$$
\begin{array}{rcc}
\sigma_{i,i-1}\colon P_i & \to & P_{i-1}\\
\cup & & \cup\\
E_i & \to & B_{i-1}
\end{array}
$$
be the resolution of the non canonical singularity of the pair
$\Box$, where $P_1=P^+$, $\sigma_{1,0}=\varphi_P$, $E_1=E_P$,
$\sigma_{2,1}$ is the blow up of the subvariety $\Lambda=B_1$, and
in general, $B_{i-1}$ is the centre of the fixed non canonical
singularity of the pair $\Box$ on $P_{i-1}$,
$E_i=\sigma^{-1}_{i,i-1}(B_{i-1})$ is the exceptional divisor,
finally, $i=1,\dots,K$ and $E_K$ realizes the fixed non canonical
singularity. Let $\Gamma$ be the oriented graph of that
resolution, that is, its set of vertices is the set of exceptional
divisors
$$
E_1,\dots,E_K,
$$
and the vertices $E_i$ and $E_j$ are joined by an oriented edge
(an {\it arrow}; notation: $i\to j$), if and only if $i>j$ and
$B_{i-1}$ is contained in the strict transform $E^{i-1}_j$ of the
exceptional divisor $E_j$ on $P_{i-1}$, see \cite{Pukh98b} or
\cite[Chapter 2]{Pukh07b}, also \cite[Chapter 2]{Pukh13a} for the
details. By the symbol $p_{ij}$ we denote the number of paths from
the vertex $E_i$ to the vertex $E_j$, if $i\neq j$; we set
$p_{ii}=1$. The fact that $E_K$ realizes a non canonical
singularity of the pair $\Box$, means that the inequality of the
Noether-Fano type holds:
\begin{equation}\label{04.06.3}
\sum^K_{i=1}p_{Ki}\nu_i>n\left(3p_{K1}+
\sum^K_{i=2}p_{Ki}\delta_i\right),
\end{equation}
where $\nu_i=\mathop{\rm mult}_{B_{i-1}}\Sigma^{i-1}$, and
$\delta_i=\mathop{\rm codim}B_{i-1}-1$ is the discrepancy of $E_i$
with respect to $P_{i-1}$. By linearity of the inequality
(\ref{04.06.3}) we may assume that $\nu_K>n$ (if $\nu_K\leq n$,
then $E_{K-1}$ is a non canonical singularity of the pair $\Box$
and $K$ can be replaced by $K-1$). Set
$$
L=\mathop{\rm max}\{2\leq i\leq K\,|\,\mathop{\rm
codim}B_{i-1}\geq3\}.
$$
The graph $\Gamma$ breaks into the {\it lower part} with the
vertices $E_1,\dots,E_L$ and the {\it upper part} with the
vertices $E_{L+1},\dots,E_K$. Now let us the well known trick of
{\it removing arrows} (see, for instance, \cite[\S 4]{Pukh10} or
\cite[Chapter 2]{Pukh13a} for the details): let us remove all
arrows that go from the vertices of the upper part to the vertex
$E_1$, if such arrows exist. This operation does not change the
numbers $p_{K2},\dots,p_{KK}$, but, generally speaking, decreases
the number of paths from $E_K$ to $E_1$. Set $p_i=p_{Ki}$ for
$i=2,\dots,K$ and let $p_1$ be the number of paths from $E_K$ to
$E_1$ in the {\it modified} graph. Since $\nu_1\leq3n$, the
inequality (\ref{04.06.3}) remains true:
\begin{equation}\label{04.06.4}
\sum^K_{i=1}p_i\nu_i>n\left(3p_1+\sum^K_{i=2}p_i\delta_i\right),
\end{equation}
in addition, the modification of the graph $\Gamma$ yields the
estimate
\begin{equation}\label{08.07.2}
p_1\leq\sum^L_{i=2}p_i.
\end{equation}
Set $Z^i_P$ to be the strict transform of the cycle $Z_P$ on
$P_i$, $i=1,\dots,L$; in particular, $Z^1_P=Z^+_P$. Set also for
$i=1,\dots,L$
$$
m_i=\mathop{\rm mult}\nolimits_{B_{i-1}}Z^{i-1}_P.
$$
Applying the technique of counting multiplicities (see, for
example \cite[Proposition 2.11]{Pukh07b} or \cite[Chapter
2]{Pukh13a}), we obtain the inequality
$$
\sum^L_{i=1}p_im_i\geq\sum^K_{i=1}p_i\nu^2_i,
$$
whence in the standard way (computing the minimum of the quadratic
form $\sum p_i\nu^2_i$ on the hyperplane, which we obtain,
replacing the inequality sign in (\ref{04.06.4}) by the equality
sign) we deduce the estimate
$$
\left(\sum^K_{i=1}p_i\right)\left(\sum^L_{i=1}p_im_i\right)>
\left(3p_1+\sum^K_{i=2}p_i\delta_i\right)^2n^2.
$$
Now set
$$
\Sigma_0=\sum_{\delta_i=3,i\geq 2}p_i,\quad\Sigma_1=
\sum_{\delta_i=2}p_i,\quad\Sigma_2=\sum_{\delta_i=1}p_i,
$$
so that, in particular, $p_1\leq\Sigma_0+\Sigma_1$. Taking into
account that the multiplicities $m_i$ do not increase, we obtain
the inequality
\begin{equation}\label{04.06.5}
(p_1+\Sigma_0+\Sigma_1+\Sigma_2)(p_1m_1+(\Sigma_0+\Sigma_1)m_2)>
(3p_1+3\Sigma_0+2\Sigma_1+\Sigma_2)^2n^2.
\end{equation}
Recall that $m_1=\mathop{\rm mult}_oZ_P$ and $m_2=\mathop{\rm
mult}_{\Lambda}Z^+_P$ are precisely the multiplicities, which we
are interested in, and we prove the inequality $m_1+m_2>12n^2$. By
linearity in $m_1,m_2$ and the last inequality (that is, the
inequality (\ref{04.06.1})) and the inequality (\ref{04.06.5}), it
is sufficient to check that the estimate (\ref{04.06.5}) does not
hold for $m_1=8n^2$, $m_2=4n^2$ and for $m_1=12n^2$, $m_2=0$.
Since $p_1\leq\Sigma_0+\Sigma_1$, it is sufficient to consider the
first case. Setting in (\ref{04.06.5}) $m_1=8n^2$ and $m_2=4n^2$,
cancelling $n^2$ and moving everything to the right hand side, we
obtain the inequality
$$
0>\Phi(p_1,\Sigma_0,\Sigma_1,\Sigma_2)
$$
where
$$
\Phi(s,t_0,t_1,t_2)=(s-t_2)^2+6st_0+5t^2_0+4t_0t_1+2t_0t_2.
$$
We obtained a contradiction, which proves the inequality
(\ref{04.06.1}).\vspace{0.1cm}

Now let us show the inequality (\ref{04.06.2}). The arguments are
completely similar to those above, with the only difference that
the coefficient at $p_1$ in the Noether-Fano inequality is 4, the
elementary discrepancies can take four, not three values, that is,
$\delta_i\in \{1,2,3,4\}$, so that there are, generally speaking,
four groups of vertices of the graph $\Gamma$ and we must set
$$
\Sigma_0=\sum_{\delta_i=4,i\geq
2}p_i,\quad\Sigma_1=\sum_{\delta_i=3}p_i,\quad\Sigma_2=
\sum_{\delta_i=2}p_i,\quad\Sigma_3=\sum_{\delta_i=1}p_i,
$$
and the inequality $p_1\leq\Sigma_0+\Sigma_1+\Sigma_2$ holds. The
technique of counting multiplicities gives the following estimate,
which is similar to the inequality (\ref{04.06.5}):
\begin{equation}\label{05.06.1}
\begin{array}{c}
(p_1+\Sigma_0+\Sigma_1+\Sigma_2+\Sigma_3)(p_1m_1+(\Sigma_0+
\Sigma_1+\Sigma_2)m_2)> \\ \\
>(4p_1+4\Sigma_0+3\Sigma_1+2\Sigma_2+\Sigma_3)^2n^2.
\end{array}
\end{equation}
Since $m_1\leq 12n^2$, to prove the inequality (\ref{04.06.2})
(which in the notations of the resolution of singularities takes
the form of the inequality $m_2>4n^2$), it is sufficient to check
that the inequality (\ref{05.06.1}) can not be true for
$m_1=12n^2$ and $m_2=4n^2$. Substituting these values into
(\ref{05.06.1}), cancelling $n^2$ and moving everything to the
right hand side, we get the inequality
$$
0>\Phi(p_1,\Sigma_0,\Sigma_1,\Sigma_2,\Sigma_3),
$$
where
$$
\Phi(s,t_0,t_1,t_2,t_3)=(2s-t_3)^2+(\dots),
$$
where in the brackets we have a quadratic form in
$s,t_0,t_1,t_2,t_3$ with nonnegative coefficients. We obtained a
contradiction, proving the inequality
(\ref{04.06.2}).\vspace{0.1cm}

The remaining inequalities of Proposition 3.7 are shown word for
word in the same way as the inequality (\ref{04.06.1}) or
(\ref{04.06.2}), depending on the type of the
inequality.\vspace{0.1cm}

Proof of Proposition 3.7 is complete.\vspace{0.1cm}

The further work, completing the proof of Theorem 5, is organized
in the following way: we exclude the cases 1.1-2.3, inspecting all
geometric possibilities.

%%%%%%%%%%%%%%%%%%%%%%%%%%%%%%%%%%%%%%%%%%%%%%%%%%%%%%%%%%%%%%%%%
%%%%%%%%%%%%%%%%%%%%%%%%%%%%%%%%%%%%%%%%%%%%%%%%%%%%%%%%%%%%%%%%%
%%%%%%%%%%%%%%%%%%%%%%%%%%%%%%%%%%%%%%%%%%%%%%%%%%%%%%%%%%%%%%%%%
%%%%%%%%%%%%%%%%%%%%%%%%%%%%%%%%%%%%%%%%%%%%%%%%%%%%%%%%%%%%%%%%%
%%%%%%%%%%%%%%%%%%%%%%%%%%%   SECTION 4

\section{Infinitely near case. II.\\ Exclusion of the linear case}

In this section we prove that the cases 1.1 and 1.2 do not
realize: it is sufficient to exclude the first one, which
immediately implies that the second one is
impossible.\vspace{0.3cm}

{\bf 4.1. Decomposition of an effective cycle.} Let us forget for
a moment about the proof of Theorem 5 and consider one very simple
construction which will be used below many times. Let $X$ be an
arbitrary algebraic variety, $Y\subset X$ an irreducible
subvariety and $Z$ an effective cycle of codimension two on $X$.
Assume first that $\mathop{\rm codim} (Y\subset X)\leq 2$, that
is, $Y$ is a prime Weil divisor on $X$ or an irreducible
subvariety of codimension two.\vspace{0.1cm}

{\bf Definition 4.1.} We say that the presentation
$$
Z=Z_0+Z_1
$$
is a $Y$-{\it decomposition} of the cycle $Z$, if both cycles
$Z_0$, $Z_1$ are effective and an irreducible component of the
cycle $Z$ is contained in $Z_0$ (respectively, in $Z_1$) if and
only if it is contained in $Y$ (respectively, not contained in
$Y$).\vspace{0.1cm}

Assume now that $\mathop{\rm codim} (Y\subset X)\geq
3$.\vspace{0.1cm}

{\bf Definition 4.2} We say that the presentation
$$
Z=Z_0+Z_1
$$
is a $Y$-{\it decomposition} of the cycle $Z$, if both cycles
$Z_0$, $Z_1$ are effective and an irreducible component of the
cycle $Z$ is contained in $Z_0$ (respectively, in $Z_1$) if and
only if it does not contain $Y$ (respectively, does contain
$Y$).\vspace{0.1cm}

Note that the definitions are not symmetric.\vspace{0.1cm}

Let us come back to the proof of Theorem 5.\vspace{0.3cm}

%%%%%%%%%%%%%%%%%%%%%%%%%%%%%%%%%%%%%%%%%%%%%%%%%%%%%%%%%%
%%%%%%%%%%%%%%%%%%%%%%%%%%%%%%%%%%%%%%%%%%%%%%%%%%%%%%%%%%
%%%%%%%%%%%%   subsection 4.2

{\bf 4.2. Restriction onto a hyperplane section.} The main result
of this section is the following\vspace{0.1cm}

{\bf Proposition 4.1.} {\it The case 1.1 does not take
place.}\vspace{0.1cm}

{\bf Proof.} Assume the converse: the case 1.1 takes place. Our
purpose is to get a contradiction. We will do it in several steps,
since the case under consideration is the hardest of the six ones.
We use both inequalities of Proposition 3.7 for the case 1.1
without special comments.\vspace{0.1cm}

First of all, let us repeat the operation of restricting onto a
hyperplane section that was used in Sec. 3.\vspace{0.1cm}

Let $R\subset P$ be a general hyperplane section, such that
:
\begin{itemize}

\item $o\in R$, the variety $R$ is non-singular at that point,

\item the hyperplane $E_R=R^+\cap E_P$ in $E_P$ contains
the subspace $\Lambda$.

\end{itemize}

Let us restrict the system $\Sigma_P$ onto $R$ and obtain a mobile
linear system $\Sigma_R$ on the hypersurface $R\subset\langle
R\rangle\cong{\mathbb P}^{M-1}$ with the self-intersection
$Z_R=Z_P|_R$, satisfying the estimates
$$
\mathop{\rm mult}\nolimits_oZ_R+ \mathop{\rm
mult}\nolimits_{\Lambda}Z^+_R>12n^2
$$
and
$$
\mathop{\rm mult}\nolimits_{\Theta}Z_R>4n^2.
$$
The advantage of this situation is that the subspaces
$\Theta\subset\Lambda\subset E_R$ are of codimension 3 and 2,
respectively. Let
$$
Z_R=Z_0+Z_1
$$
be the $T_R$-decomposition of the cycle $Z_R$, where
$T_R=(T_oR)\cap R$ is the tangent hyperplane section at the point
$o$. Set
$$
d_0=\frac{1}{Mn^2}\mathop{\rm deg}Z_0,\quad
d_1=\frac{1}{Mn^2}\mathop{\rm deg}Z_1,
$$
$$
\mu_0=\frac{1}{n^2}\mathop{\rm mult}\nolimits_oZ_0,\quad
\mu_1=\frac{1}{n^2}\mathop{\rm mult}\nolimits_oZ_1.
$$
We obtain the equality
\begin{equation}\label{25.06.1}
d_0+d_1=4
\end{equation}
and the inequality
\begin{equation}\label{25.06.2}
\mu_0+\mu_1>8.
\end{equation}
Furthermore, set $\lambda_1=\frac{1}{n^2}\mathop{\rm
mult}_{\Lambda}Z^+_1$, where $Z^+_1$ is the strict transform of
the cycle $Z_1$ on $R^+$, so that the following inequality holds:
\begin{equation}\label{25.06.3}
\mu_0+\mu_1+\lambda_1>12.
\end{equation}
Proposition 3.3 implies that the multiplicities $\mu_i$ can be
estimated in terms of the degrees $d_i$ in the following way: for
$M\geq 18$ the inequality
\begin{equation}\label{25.06.4}
\mu_0\leq 3d_0
\end{equation}
holds, for $M\leq 17$ a weaker estimate is true:
\begin{equation}\label{25.06.5}
\mu_0\leq\frac{8M}{3(M-2)}d_0.
\end{equation}
Since none of the components of the cycle $Z_1$ is contained in
the tangent section $T_R=\{q_1|_R=0\}=R\cap T_oR$, by the part
(ii) of Proposition 3.3 the inequality
\begin{equation}\label{25.06.6}
\mu_1\leq\frac{2M}{M-2}d_1
\end{equation}
holds. Finally, it is obvious that the following estimate holds:
\begin{equation}\label{25.06.7}
\mu_1\geq\lambda_1.
\end{equation}

The system of six equations and inequalities
(\ref{25.06.1}-\ref{25.06.7}) (it is six, because depending on
whether $M\geq 18$ or $M\leq 17$, we choose the inequality
(\ref{25.06.4}) or (\ref{25.06.5})) forms the {\it first} system
of relations for the five parameters introduced
above.\vspace{0.1cm}

Since $q_2|_{\Theta}\not\equiv 0$, the components of the cycle
$Z_R$, the strict transforms of which contain the linear subspace
$\Theta$, can not be contained in $T_R$. For that reason,
$\mathop{\rm mult}_o Z_1\geq \mathop{\rm mult}_{\Theta} Z^+_1
>4n^2$, so that the following inequality holds:
\begin{equation}\label{02.10.2013.1}
\mu_1 > 4.
\end{equation}
\vspace{0.1cm}

%%%%%%%%%%%%%%%%%%%%%%%%%%%%%%%%%%%%%%%%%%%%%%%%%%%%%%%%%%%%%%%
%%%%%%%%%%%%%%%%%%%%%%%%%%%%%%%%%%%%%%%%%%%%%%%%%%%%%%%%%%%%%%%
%%%%%%%%%%%%%%%%%%   subsection 4.3

{\bf 4.3. Additional estimates for the cycle $Z_1$.} Now let us
consider the cycle $Z_1$, the most important part of the
self-intersection $Z_R$, since it contains the linear subspace
$\Lambda$. First of all, none of the components of the cycle $Z_1$
is contained in the tangent section $T_R=R\cap T_o R$ and for that
reason $(Z_1\circ T_R)$ is an effective cycle of codimension 2 on
the hypersurface $T_R$. The latter has a quadratic singularity at
the point $o$, so that
$$
\{q_2|_{T_{o}R}=0\}
$$
is its tangent cone at that point. Its projectivization will be
denoted by the symbol $Q_R$. Obviously,
$$
Q_R=T^+_R\cap E.
$$
By the condition (R2) the intersection $[Q_R\cap\Lambda]$ is an
irreducible quadric; it is subvariety of codimension two on $Q_R$.
Now let us compute the multiplicity $\mathop{\rm mult}_o(Z_1\circ
T_R)$. By the rules of the intersection theory (see \cite{Ful} or
\cite[Chapter 2]{Pukh13a}), write
$$
(Z^+_1\circ T^+_R)=(Z_1\circ T_R)^++Z_Q,
$$
where $Z_Q$ is an effective divisor on the quadric $Q_R$ (outside
$Q_R$ the effective cycles $(Z^+_1\circ T^+_R)$ and $(Z_1\circ
T_R)^+$ obviously coincide). Now we have
$$
\mathop{\rm mult}\nolimits_o(Z_1\circ T_R)=2\mu_1n^2+\mathop{\rm
deg}Z_Q.
$$
Setting $\mu_{2}=\frac{1}{n^2}\mathop{\rm mult}_o(Z_{1}\circ
T_R)$, and $\mathop{\rm deg}Z_Q=2\delta n^2$, we obtain the
equality
\begin{equation}\label{27.06.1}
\mu_{2}=2\mu_1+2\delta.
\end{equation}
Obviously, $\mathop{\rm deg} (Z_1\circ T_R)=\mathop{\rm deg} Z_1$.
Furthermore, set
$$
\lambda_2=\frac{1}{n^2}\mathop{\rm mult}\nolimits_{[Q_R\cap
\Lambda]} (Z_1\circ T_R)^+.
$$
Since the following inequality is obviously true:
$$
\mathop{\rm mult}\nolimits_{[Q_R\cap\Lambda]}(Z^+_1\circ
T^+_R)\geq \mathop{\rm mult}\nolimits_{\Lambda}Z^+_1,
$$
by Corollary 4.1, which is shown below, we get the estimate
\begin{equation}\label{27.06.4}
\lambda_2\geq\lambda_1-\delta.
\end{equation}
Finally, since $\mathop{\rm mult}\nolimits_o(Z_1\circ T_R) \geq
\mathop{\rm mult}\nolimits_{[Q_R\cap\Lambda]}(Z_1\circ T_R)^+$, we
obtain the estimate
\begin{equation}\label{27.06.5}
\mu_{2}\geq\lambda_2.
\end{equation}
Now we have to take into account the input of the infinitely near
subvariety $[Q_R\cap\Lambda]$.\vspace{0.1cm}

{\bf Proposition 4.2.} {\it The following estimate holds:}
\begin{equation}\label{27.06.6}
\mu_{2}+2\lambda_2\leq 4\frac{M}{M-3}d_{1}.
\end{equation}

{\bf Proof.} Let $H_R$ be the class of a hyperplane section of the
hypersurface $R\subset \langle R\rangle\cong{\mathbb P}^{M-1}$.
Consider the pencil $|H_R-\Lambda|$ of hyperplane sections of $R$,
defined by the condition: for $S\in |H_R-\Lambda|$ we have
$S^+\supset\Lambda$. The base set $\Delta_R$ of the pencil
$|H_R-\Lambda|$ is an irreducible subvariety of codimension two in
$T_R$, of codimension 3 in $R$; more precisely,
$\Delta_R\subset\langle\Delta_R\rangle\cong{\mathbb P}^{M-3}$ is a
hypersurface of degree $M$, where the linear span
$\langle\Delta_R\rangle$ is determined by the condition
$$
\langle\Delta_R\rangle^+\cap E=\Lambda.
$$

Now let $Y$ be an arbitrary subvariety of codimension 2 in $R$.
For a general divisor $S\in|H_R-\Delta|$ we have $Y\not\subset S$,
so that $(Y\circ S)$ is an effective cycle of codimension  3 on
$R$. By construction,
$$
\mathop{\rm mult}\nolimits_o(Y\circ S)\geq \mathop{\rm
mult}\nolimits_oY+2\mathop{\rm
mult}\nolimits_{[Q_R\cap\Lambda]}Y^+
$$
(since $\mathop{\rm deg}[Q_R\cap\Lambda]=2$). However, $S$ is a
section of the singular hypersurface $V\cap T_oV$ by a linear
subspace of codimension 3, so that, applying the regularity
condition (R1) and arguing in word for word the same way as in the
proof of Proposition 3.3, by means of the technique of
hypertangent divisors, applied to the cycle $(Y\circ S)$, we
obtain the estimate
\begin{equation}\label{27.06.7}
\mathop{\rm mult}\nolimits_oY+ 2\mathop{\rm
mult}\nolimits_{[Q_R\cap\Lambda]}Y^+\leq \frac{4}{M-3}\mathop{\rm
deg}Y
\end{equation}
(recall that on $S$ the cycle $(Y\circ S)$ has codimension 2, so
that this cycle can be considered as an effective cycle of
codimension 3 on a section of the hypersurface $V$ by a linear
subspace of codimension 3, which by the condition (R1) satisfies
the regularity condition).\vspace{0.1cm}

The inequality (\ref{27.06.6}) follows from (\ref{27.06.7}) in a
trivial way.\vspace{0.1cm}

Proof of Proposition 4.2 is complete.\vspace{0.3cm}

%%%%%%%%%%%%%%%%%%%%%%%%%%%%%%%%%%%%%%%%%%%%%%%%%%%%%%%%%%%%%%%%%%%
%%%%%%%%%%%%%%%%%%%%%%%%%%%%%%%%%%%%%%%%%%%%%%%%%%%%%%%%%%%%%%%%%%%
%%%%%%%%%%%%%%%%%%%%   subsection 4.4

{\bf 4.4. On the multiplicities of subvarieties on a quadric.} Let
us put off the proof of Theorem 5 and show the fact about
multiplicities of subvarieties on a quadric hypersurface that was
used in Subsection 4.3. Let $Q\subset{\mathbb P}^N$ be an
irreducible quadric, $\dim \mathop{\rm Sing} Q=s_Q$, and $Y\subset
X\subset Q$ irreducible subvarieties.\vspace{0.1cm}

{\bf Proposition 4.3.} {\it Assume that the inequality
\begin{equation}\label{24.06.2}
\dim X+ \dim Y> N+s_Q+1
\end{equation}
holds. The then following estimate is true:}
\begin{equation}\label{24.06.1}
\mathop{\rm mult}\nolimits_Y X\leq \frac12 \mathop{\rm deg} X.
\end{equation}

{\bf Proof.} Assume the converse:
$$
2\mathop{\rm mult}\nolimits_Y X >\mathop{\rm deg} X.
$$
By the assumption on the dimensions $Y\not\subset \mathop{\rm
Sing} Q$. Take an arbitrary point $p\in Y\setminus \mathop{\rm
Sing} Q$.\vspace{0.1cm}

{\bf Lemma 4.1.} {\it The variety $X$ is contained in the tangent
hyperplane} $T_pQ$.\vspace{0.1cm}

{\bf Proof.} Assume the converse: $X\not\subset T_pQ$. Then the
effective cycle $(X\circ (T_pQ\cap Q))$ is well defined. Its
degree is equal to $\mathop{\rm deg} X$ and its multiplicity at
the point $p$ is at least $2\mathop{\rm mult}\nolimits_p
X>\mathop{\rm deg} X$, which is impossible. Q.E.D. for the
lemma.\vspace{0.1cm}

Therefore, the following inclusion takes place
$$
X\subset\mathop{\bigcap}\limits_{p\in Y\setminus \mathop{\rm Sing}
Q} T_pQ.
$$
Proof of the proposition will be complete, if we show that the
dimension of the right hand side of the inclusion is strictly
smaller than $\dim X$. Note that $\mathop{\rm Sing} Q\subset
{\mathbb P}^N$ is a linear subspace and for any non-singular point
$p\in Q$ we have $\mathop{\rm Sing} Q\subset T_p Q$.\vspace{0.1cm}

Consider the section $Q^*$ of the quadric $Q$ by a general linear
subspace of codimension $s_Q+1$ (in particular, not meeting
$\mathop{\rm Sing} Q$). The quadric $Q^*$ is non-singular. Let
$Y^*\subset X^*$ be the corresponding sections of the varieties
$Y$ and $X$.\vspace{0.1cm}

Obviously, $Y^*$ contains at least
$$
\dim Y^* =\dim Y-s_Q -1
$$
linearly independent points, so that the linear space
$$
\mathop{\bigcap}\limits_{p\in Y^*} T_p Q^*
$$
has the dimension not higher than the number
$$
N-(s_Q+1)-(\dim Y-s_Q-1)=N-\dim Y.
$$
Therefore, $\dim X^*=\dim X-s_Q-1\leq N-\dim Y$. However, by
assumption the opposite inequality holds. This contradiction
completes the proof of Proposition 4.3.\vspace{0.1cm}

{\bf Corollary 4.1.} {\it Let $o\in V$ be an arbitrary point,
$\Pi\subset T_oV$ a linear subspace of codimension two in the
vector tangent space $T_oV\cong {\mathbb C}^M$, ${\mathbb
P}(\Pi)\cong {\mathbb P}^{M-3}$ its projectivization. Let $X$, $Y$
be irreducible subvarieties of codimension 1 and 2 on the quadric
hypersurface $Q=\{q_2|_{{\mathbb P}(\Pi)}=0\}\subset {\mathbb
P}(\Pi)$. Then the estimate (\ref{24.06.1}) holds.}\vspace{0.1cm}

{\bf Proof.} By the regularity condition (R2) for the quadric $Q$
we have the estimate $s_Q\leq [\sqrt{M}]+1$. Now by easy
computations we see that the inequality (\ref{24.06.2}) holds.
Applying Proposition 4.3, we complete the proof.\vspace{0.3cm}

%%%%%%%%%%%%%%%%%%%%%%%%%%%%%%%%%%%%%%%%%%%%%%%%%%%%%%%%%%%
%%%%%%%%%%%%%%%%%%%%%%%%%%%%%%%%%%%%%%%%%%%%%%%%%%%%%%%%%%%
%%%%%%%%%%%%%%%%%%%   subsection 4.5

{\bf 4.5. Exclusion of the linear case.} Let us complete the proof
of Proposition 4.1. The six linear equations and inequalities
(\ref{27.06.1}-\ref{27.06.6}) form the {\it second} system of
relations for now 5+3=8 parameters $d_*$, $\mu_*$, $\lambda_*$ and
$\delta$. Joining the first and second systems of relations,
adding to them the inequality (\ref{02.10.2013.1}), we obtain 11
linear equations and inequalities for 8 nonnegative real
parameters. replacing the strict inequalities everywhere by the
non-strict ones, we obtain 11+8=19 linear equations and non-strict
linear inequalities, defining some (obviously, compact) convex
subset
$$
\Xi\subset{\mathbb R}^{8}.
$$
%\vspace{0.1cm}

{\bf Proposition 4.4.} {\it The set $\Xi$ is empty.}\vspace{0.1cm}

{\bf Proof.} It is sufficient to apply any computer program to
solve a suitable problem of linear programming, for example
$$
\mu_0\to\mathop{\rm max}_{\Xi}.
$$
For MAPLE the corresponding command can be written in the
following way:\vspace{0.1cm}

\noindent {\tt >with(Optimization):}

\noindent {\tt
>M:=15:LPSolve(m0,\{m0+m1>=8,d0+d1=4,m0<=(8*M/(3*(M-2)))*d0,}

\noindent {\tt
m1<=2*(M/(M-2))*d1,m0+m1+l1>=12,m1>=l1,m1>=4,m0>=0,}

\noindent {\tt
d0>=0,d1>=0,l1>=0,m2=2*m1+2*de,l2>=l1-de,m2+2*l2<=4*(M/(M-3))*d1,}

\noindent {\tt m2>=l2,de>=0,l2>=0\},
maximize);}\vspace{0.1cm}

and its application gives the following result:\vspace{0.1cm}

\noindent {\tt Error, (in Optimization:-LPSolve) no feasible
solution found}\vspace{0.1cm}

This completes the proof of Proposition 4.4. (The convex set $\Xi$
is defined by 17 linear inequalities in the affine subspace of
codimension two
$$
\{d_0+d_1=4,\,\, \mu_2=2\mu_1+2\delta\}\subset {\mathbb R}^8,
$$
so that to solve the problem of linear programming, one needs to
inspect a finite set of points of bounded cardinality. Each of
these points is checked for being a point of the set $\Xi$.
Therefore, the proof of Proposition 4.4 can be given over to the
computer.) \vspace{0.1cm}

Therefore, the first and second systems of relations, obtained
above, define the empty set in ${\mathbb R}^{8}$. Therefore, the
case 1.1 does not take place. Proposition 4.1 is shown. Q.E.D.
\vspace{0.1cm}

{\bf Corollary 4.2.} {\it The case 1.2 does not take
place.}\vspace{0.1cm}

{\bf Proof.} Since a non log canonical singularity of the pair
$\Box^*$ is automatically a non log canonical singularity of the
pair $\Box$, the case 1.2 is a version of the case 1.1 (for
$\Theta$ one can take any hyperplane in $\Lambda$). Q.E.D. for the
corollary.

%%%%%%%%%%%%%%%%%%%%%%%%%%%%%%%%%%%%%%%%%%%%%%%%%%%%%%%%%%%%%%%%%%%
%%%%%%%%%%%%%%%%%%%%%%%%%%%%%%%%%%%%%%%%%%%%%%%%%%%%%%%%%%%%%%%%%%%
%%%%%%%%%%%%%%%%%%%%%%%%%%%%%%%%%%%%%%%%%%%%%%%%%%%%%%%%%%%%%%%%%%%
%%%%%%%%%%%%%%%%%%%%%%%%%%%%%%%%%%%%%%%%%%%%%%%%%%%%%%%%%%%%%%%%%%%
%%%%%%%%%%%%%%%%%%%%%%   SECTION 5

\section{Infinitely near case. III.\\ Exclusion of the non-linear case}

In this section we exclude the case 2, which completes the
exclusion of the infinitely near case (and so the proof of Theorem
5).\vspace{0.3cm}

{\bf 5.1. The case 2.1, $B$ is not contained in a quadric.} Let us
consider first the 2.1 and assume that the subvariety $B$ is not
contained in any quadric hypersurface in $E_P$. (Note, that in the
case 2.1 the subvariety $B$ is certainly not contained in the
quadric $Q_P$, since $B\supset \Theta$ and $\Theta\not\subset
Q_P$.)\vspace{0.1cm}

{\bf Proposition 5.1} {\it The following inequality holds:}
$$
5\mathop{\rm mult}\nolimits_B\Sigma^+_P\leq 2\mathop{\rm
mult}\nolimits_o\Sigma_P.
$$

{\bf Proof.} Recall that the multiplicity $\mathop{\rm
mult}_o\Sigma_P$ is denoted by the letter $\nu$ and set
$\nu_B=\mathop{\rm mult}_B\Sigma^+_P$. In terms of the resolution
of the maximal singularity of the system $\Sigma_P$ we have
$\nu=\nu_1$ and $\nu_B=\nu_2$. Assume that the opposite inequality
holds: $5\nu_B>2\nu$. Let us show that this implies that the
linear system $\Sigma^+_P$ can not be mobile, more precisely, the
following inclusion takes place:
$$
E_P\subset\mathop{\rm Bs}\Sigma^+_P.
$$
This contradiction implies the claim of our
proposition.\vspace{0.1cm}

Since the subvariety $B$ is not contained in any quadric
hypersurface, its degree $\mathop{\rm deg}B$ (as a subvariety of
the projective space $E_P$) is at least 5.\vspace{0.1cm}

Indeed, $d_B=\mathop{\rm deg}B\geq 3$. Furthermore, $B$ is not a
cone over a curve: otherwise, $B$ contains a linear subspace of
codimension 3 in $E_P$, which is excluded by the proof of
Proposition 4.1. Now, projecting from a point of general position
$p\in B$, we exclude the option $d_B=3$. If $d_B=4$ and $B$ has at
least one singular point of multiplicity 2 or 3, then $B$ is
contained in a hyperplane or an irreducible quadric, contrary to
the assumption. Since a non-singular projective subvariety of
codimension 2 and degree 4 in ${\mathbb P}^k$, $k\geq 4$, is a
complete intersection of two quadrics (this is a well known fact;
see also \cite{Holme89}), then $B$ is a complete intersection of
two quadrics, either, if $B$ is non-singular or is a cone over a
subvariety of degree 4 and dimension $\geq 2$. We have inspected
all options. Therefore, $d_B\geq 5$.\vspace{0.1cm}

Let $\Pi\subset E_P$ be a 2-plane of general position, so that
$B_{\Pi}=B\cap\Pi$ is a finite set, consisting of $d_B\geq 5$
distinct points. Let $R_1,\dots,R_m$ be all irreducible
hypersurfaces in $E_P$, containing $B$ and contained in
$\mathop{\rm Bs}\Sigma_P$, if there are any. Then $\mathop{\rm
deg}R_i\geq 3$ and the irreducible curves $R_i\cap\Pi$ are all
irreducible curves in the plane $\Pi$, contained in $\mathop{\rm
Bs}\Sigma_P$ and containing {\it at least one} point of the finite
set $B_{\Pi}$.\vspace{0.1cm}

{\bf Lemma 5.1.} {\it Neither three points of the set $B_{\Pi}$
are collinear.}\vspace{0.1cm}

{\bf Proof.} Assume the converse: there are three distinct points
$p_1,p_2,p_3\in B_{\Pi}$, lying on the line $L$. Since $\nu_B>n$
and $\nu\leq 3n$, we obtain, that $L\subset\mathop{\rm
Bs}\Sigma_P$. As we noted above, this is impossible. Q.E.D. for
the lemma.\vspace{0.1cm}

Now let us consider any 5 distinct points $p_1,\dots,p_5\in
B_{\Pi}$ and the unique conic $C\subset\Pi$, containing those
points. The restriction $\Sigma_C=\Sigma^+_P|_C$ is a linear
series of degree $2\nu$ with 5 base points of multiplicity
$\nu_B$. Since $5\nu_B>2\nu$, we have $C\subset\mathop{\rm
Bs}\Sigma_P$, which is impossible. Proof of Proposition 10.09.1 is
complete.\vspace{0.1cm}

Now we can apply the technique of counting multiplicities and
estimate the multiplicity of the self-intersection $Z_P$ at the
point $o$ and its strict transform $Z^+_P$ along the subvariety
$B$.\vspace{0.1cm}

Set $\mu=\mathop{\rm mult}_oZ_P$, $\mu_B=\mathop{\rm
mult}_BZ^+_P$.\vspace{0.1cm}

{\bf Proposition 5.2.} {\it The following inequality holds:}
\begin{equation}\label{10.09.1}
\mu+\mu_B>\frac{81}{5}n^2.
\end{equation}

{\bf Proof.} As in the proof of Proposition 3.7, fix a maximal
singularity, the centre of which on $P^+$ is a subvariety $B$ and
take its resolution. We use the standard notations, associated
with the resolution. The graph $\Gamma$ is assumed to be modified,
so that the inequality
$$
p_1\leq\Sigma_0=\sum^L_{i=2}p_i
$$
holds. We have the inequality of Noether-Fano type
\begin{equation}\label{10.09.2}
\sum^K_{i=1}p_i\nu_i>(3p_1+2\Sigma_0+\Sigma_1)n,
\end{equation}
where $\Sigma_1=\sum\limits^K_{i=L+1}p_i$, and, besides, we know
that $\nu_1\leq 3n$ and $5\nu_2\leq 2\nu_1$; the multiplicities
$\nu_i$ do not increase,
$$
\nu_2\geq\nu_3\geq\dots\geq\nu_K.
$$
By the technique of counting multiplicities, taking into account
the inequalities
$$
\mathop{\rm mult}\nolimits_{B_i}Z^i_P\geq \mathop{\rm
mult}\nolimits_{B_{i+1}}Z^{i+1}_P,
$$
we obtain the estimate
$$
p_1\mu+\Sigma_0\mu_B\geq\sum^K_{i=1}p_i\nu^2_i.
$$
For $\nu_1=\nu$ fixed, the minimum of the right hand side of the
latter inequality on the hyperplane
$$
\sum^K_{i=1}p_i\nu_i=(3p_1+2\Sigma_0+\Sigma_1)n
$$
is attained at $\nu_2=\dots=\nu_K=\theta$, where the value
$\theta$ is computed from the equation
\begin{equation}\label{10.09.3}
p_1\nu+(\Sigma_0+\Sigma_1)\theta=(3p_1+2\Sigma_0+\Sigma_1)n.
\end{equation}
Therefore, the inequality
\begin{equation}\label{12.09.2}
p_1\mu+\Sigma_0\mu_B>p_1\nu^2+(\Sigma_0+\Sigma_1)\theta^2
\end{equation}
holds. On the other hand, the equality (\ref{10.09.3}) can be
re-written in the following way:
$$
\Sigma_1=\frac{3n-\nu}{\theta-n}p_1+\frac{2n-\theta}{\theta-n}\Sigma_0.
$$
Recall that $\nu$ and $\theta$ are connected by the inequality
$5\theta\leq 2\nu$. As a result, we obtain that the sum
$\mu+\mu_B$ is strictly higher than the minimum of the function
$x+y$ on the interval, cut out by the inequalities
$$
x\geq\nu^2,\,\,x\geq y,\,\,y\geq 0
$$
on the line
$$
\{p_1x+\Sigma_0y=\Psi(\nu,\theta)\}\subset{\mathbb R}^2_{x,y},
$$
where
$$
\Psi(\nu,\theta)=p_1\left(\nu^2+\frac{3n-\nu}{\theta-n}\theta^2+
\frac{n\theta^2}{\theta-n}\Sigma_0\right).
$$
The more so, this minimum is strictly higher than the number
\begin{equation}\label{13.09.1}
\nu^2+\frac{n\theta^2}{\theta-n}.
\end{equation}
It is easy to check that the minimum of the function
(\ref{13.09.1}) on the triangle
$$
\{\theta>n,\,\,\nu\leq 3n,\,\,5\theta\leq 2\nu\}\subset{\mathbb
R}^2_{\nu,\theta}
$$
is attained for $\nu=3n$, $\theta=\frac{6}{5}n$ and is equal to
$\frac{81}{5}n^2$. This completes the proof of Proposition
5.2.\vspace{0.1cm}

The inequality (\ref{10.09.1}) is so strong that it makes it
possible to easily complete the exclusion of the case 2.1 (under
the assumption that $B$ is not contained in any quadric
hypersurface in $E_P$). Indeed, since $d_B\geq 5$, we have the
inequality
$$
\mu\geq 5\mu_B.
$$
It is easy to check that it is incompatible with the inequalities
(\ref{10.09.1}) and $\mu\leq 12n^2$. This excludes the case under
consideration (that is, the case 2.1 under the assumption that $B$
is not contained in any quadric in $E_P$).\vspace{0.3cm}

%%%%%%%%%%%%%%%%%%%%%%%%%%%%%%%%%%%%%%%%%%%%%%%%%%%%%%%%%%%%%%%%%%%
%%%%%%%%%%%%%%%%%%%%%%%%%%%%%%%%%%%%%%%%%%%%%%%%%%%%%%%%%%%%%%%%%%%
%%%%%%%%%%%%%%%%%%%%%%%   subsection 5.2

{\bf 5.2. Case 2.1, $B$ is contained in a quadric, but not in a
hyperplane.} Now let us consider the case 2.1 under the assumption
that $B$ is contained in some quadric in $E_P$, but $\langle
B\rangle=E_P$, that is, $B$ is not contained in any hyperplane in
$E_P$.\vspace{0.1cm}

{\bf Proposition 5.3.} {\it The following inequality holds:}
$$
2\mathop{\rm mult}\nolimits_B\Sigma^+_P\leq \mathop{\rm
mult}\nolimits_o\Sigma_P.
$$

{\bf Proof.} Again we write $\nu_B=\mathop{\rm mult}_B\Sigma^+_P$
and $\nu=\mathop{\rm mult}_o\Sigma_P$. Since $B$ is not contained
in a hyperplane, $\mathop{\rm Sec}(B)=E_P$. Let $L$ be a general
secant line of the variety $B$. Since the system $\Sigma^+_P$ has
no fixed components, for a general divisor $D\in\Sigma_P$ we have
$L\not\subset D^+$. Therefore,
$$
2\nu_B\leq\sum_{x\in L\cap B}(L\cdot D^+)_x\leq(L\cdot D^+)=\nu,
$$
as we claimed. The proposition is shown.\vspace{0.1cm}

{\bf Corollary 5.1.} {\it The following estimate is true:}
$\nu_B\leq\frac32n$.\vspace{0.1cm}

The following claim is an analog of Proposition 5.2 in the
situation under consideration\vspace{0.1cm}

{\bf Proposition 5.4.} {\it The following inequality holds:}
\begin{equation}\label{12.09.1}
\mu+\mu_B>(10+2\sqrt{2})n^2.
\end{equation}

{\bf Proof} is completely similar to the proof of Proposition 5.2
given above: we argue in word for word the same way and, recalling
that $\mu>8n^2$, we get that the value $\mu+\mu_B$ is strictly
higher than the minimum of the function
$$
\mathop{\rm max}(\nu^2,8n^2)+\frac{n\theta^2}{\theta-n}
$$
on the triangle
$$
\{\theta>n,\,\,\nu\leq 3n,\,\, 2\theta\leq\nu\}\subset {\mathbb
R}^2_{\nu,\theta}.
$$
This minimum is attained for $\nu=2\sqrt{2}n$, $\theta=\sqrt{2}n$
and is equal to $(10+2\sqrt{2})n^2$, which is what we need.
Q.E.D.\vspace{0.1cm}

{\bf Remark 5.1.} Since $10+2\sqrt{2}\thickapprox 12.8$, the
inequality (\ref{12.09.1}) is considerably sharper than the
aprioric inequality (\ref{29.06.1}).\vspace{0.1cm}

The estimate (\ref{12.09.1}) is essentially weaker than
(\ref{10.09.1}), however, this is compensated by the additional
geometric information about the subvariety $B$: we know that
$B\subset Q^*$, where $Q^*\neq Q_P$ is some irreducible quadric,
and moreover by assumption $B$ is not a hyperplane section of the
quadric $Q^*$.\vspace{0.1cm}

{\bf Lemma 5.2.} {\it The degree of the subvariety $B$ is at least
4.}\vspace{0.1cm}

{\bf Proof.} We must exclude the option $d_B=\mathop{\rm deg}B=3$.
Assume that this is the case. Then the rank of the quadratic form,
defining $Q^*$, is equal to 3 or 4, so that $B$ is swept out by a
one-dimensional family of linear subspaces of codimension 3 in
$E_P$. Proposition 4.1 excludes this situation. Q.E.D. for the
lemma.\vspace{0.1cm}

{\bf Corollary 5.2.} {\it The following inequality holds:}
$$
\mu\geq 4\mu_B.
$$

Now we exclude the case under consideration in the same way as we
used to exclude the case 1.1, with some simplifications. Let
$Z_P=Z_0+Z_1$ be the $T_P$-decomposition of the cycle $Z_P$. Since
$B\not\subset Q_P=T^+_P\cap E_P$, we have $\mathop{\rm
mult}_BZ^+_1=\mu_B$. Now, introducing the normalized parameters
$d_i$, $\mu_i$, $i=0,1$, and $\lambda_1$, we obtain for them the
system of the following inequalities: (\ref{25.06.1}),
(\ref{25.06.2}), instead of (\ref{25.06.4}) and (\ref{25.06.5}) we
have the estimate
$$
\mu_0\leq\max\left( 3,\frac{8M}{3(M-1)}\right)\, d_0,
$$
instead of (\ref{25.06.6}) we have the estimate
$$
\mu_1\leq \frac{2M}{M-1}\, d_1,
$$
finally, instead of (\ref{25.06.3}) we have the estimate
$$
\mu_0+\mu_1+\lambda_1 > 10+2\sqrt{2}
$$
and instead of (\ref{25.06.7}) the stronger estimate
$$
\mu_1\geq 4\lambda_1.
$$
Using MAPLE, it is easy to check that this system of linear
equations and inequalities has no solutions already for $M\geq 5$.
This completes the exclusion of the case 2.1 under the assumption
that $\langle B\rangle =E_P$.\vspace{0.3cm}

%%%%%%%%%%%%%%%%%%%%%%%%%%%%%%%%%%%%%%%%%%%%%%%%%%%%%%%%%%%%%%%%%%
%%%%%%%%%%%%%%%%%%%%%%%%%%%%%%%%%%%%%%%%%%%%%%%%%%%%%%%%%%%%%%%%%%
%%%%%%%%%%%%%%%%%%%   subsection 5.3

{\bf 5.3. The case 2.1, $B$ is contained in a hyperplane.} Assume
that $B$ is contained in some hyperplane $\Pi\subset E_P$. By
Proposition 4.1, $B$ is a hypersurface of degree $d_B\geq 2$ in
$\Pi$. Consider the linear system $|H_P-\Pi|$, that is, the pencil
of hyperplane sections, the base set of which is the intersection
$\Delta$ of the tangent section $T_P$ with the hyperplane in
$\langle P\rangle$ that has $\Pi$ as the tangent cone. Let
$Z_P=Z_0+Z_1$ be such a decomposition of the cycle $Z_P$, that
$Z^+_P=Z^+_0+Z^+_1$ is the $B$-decomposition of the effective
cycle $Z^+_P$.\vspace{0.1cm}

Let $R\in|H_P-\Pi|$ be a general divisor. For the effective cycle
$(Z_1\circ R)$ we have:
$$
\begin{array}{c}
\mathop{\rm deg}(Z_1\circ R)=\mathop{\rm deg}Z_1,\\ \\ \mathop{\rm
mult}\nolimits_o(Z_1\circ R)\geq
\mathop{\rm mult}\nolimits_oZ_1+2\mathop{\rm mult}\nolimits_BZ^+_1,\\
\end{array}
$$
since $d_B\geq 2$. However, in the case under consideration
$(Z_1\circ R)$ is an effective cycle of codimension two on the
hyperplane section $R$, which itself satisfies the regularity
conditions, and for that reason the inequality
$$
\frac{\mathop{\rm mult}\nolimits_o}{\mathop{\rm deg}}(Z_1\circ R)
\leq\mathop{\rm max}\left(\frac{3}{M},\frac{8}{3(M-2)}\right)
$$
holds; the right hand side for $M\geq 18$ does not exceed $3/M$.
Taking into account that $\mathop{\rm mult}_oZ_0\leq
\frac{3}{M}\mathop{\rm deg}Z_0$, we obtain a contradiction with
the aprioric inequality (\ref{29.06.1}). This excludes the case
under consideration for $M\geq 18$.\vspace{0.1cm}

{\bf Remark 5.3.} In the argument given above we used the fact
that $B\not\subset Q_P$: it is for that reason that the
scheme-theoretic intersection $(Z_1\circ R)$ is well defined.
However, if $B\subset Q_P$, then $B=\Pi\cap Q_P$. Set
$\Delta=\mathop{\rm Bs}|H_P-\Pi|$ (see above). Obviously,
$$
\mathop{\rm deg}\Delta=M,\quad\mathop{\rm mult}\nolimits_o
\Delta=2
$$
(because $\Delta^+\cap E_P=B$), so that writing
$$
где Z_1=a\Delta+Z_*,
$$
where $a\in{\mathbb Z}_+$ and $Z_*$ does not contain $\Delta$ as a
component, we repeat the previous argument and come to a
contradiction for $M\geq 18$.\vspace{0.1cm}

If we use all the information available, we can exclude the case
under consideration for smaller values of $M$ as well. Namely,
write $Z_P=Z_0+Z_1$, where $Z^+_P=Z^+_0+Z^+_1$ is the
$\Theta$-decomposition of the effective cycle $Z^+_P$. The cycle
$(Z_1\circ R)$ is well defined for a general divisor
$R\in|H_P-\Pi|$. Furthermore, write
$$
(Z_1\circ R) =Z_{10}+Z_{11},
$$
where the strict transform of this equality on $P^+$ is the
$\Theta$-decomposition of the cycle $(Z_1\circ R)^+$. The cycle
$Z_{10}$ satisfies the estimate
$$
\frac{\mathop{\rm mult}\nolimits_o}{\mathop{\rm deg}}Z_{10}\leq
\frac{8}{3(M-2)}d_{10},
$$
but for $Z_{11}$ a much stronger inequality holds:
$$
\frac{\mathop{\rm mult}\nolimits_o}{\mathop{\rm deg}}Z_{11}\leq
\frac{2}{M-2}d_{11},
$$
since none of the components of the cycle $Z_{11}$ is contained in
the tangent section $T_R=R\cap T_oR$, so that we can form the
effective cycle $(Z_{11}\circ T_R)$ and then apply to this cycle
of codimension two on $T_R$ the technique of hypertangent
divisors. Finally, setting
$$
\xi_1=\frac{1}{n^2}\mathop{\rm mult_{\Theta}}Z^+_P,\quad\xi_2=
\frac{1}{n^2}\mathop{\rm mult_{\Theta}}Z^+_{11},
$$
we get the following system of linear equations and inequalities:
(\ref{25.06.1}), (\ref{25.06.2}), (\ref{25.06.3}),
(\ref{25.06.4}),  and also
$$
\begin{array}{lll}
\mu_{10}+\mu_{11}=\mu_1+2\lambda_1+\delta_1, & \xi_1>4, & \\ \\
d_{10}+d_{11}=d_1, & \xi_2\geq\xi_1-\delta_1, &
\mu_{11}\geq\xi_2,\\ \\
\displaystyle \mu_{10}\leq\frac{8M}{3(M-2)}d_{10}, & \displaystyle
\mu_{11}\leq\frac{2M}{(M-2)}d_{11}. &
\end{array}
$$
Applying MAPLE we see that this system is incompatible (even when
we replace all strict inequalities by the non-strict ones) already
for $M\geq 11$. This completes the exclusion of the case
2.1.\vspace{0.3cm}

%%%%%%%%%%%%%%%%%%%%%%%%%%%%%%%%%%%%%%%%%%%%%%%%%%%%%%%%%%%%%%%%%%%%
%%%%%%%%%%%%%%%%%%%%%%%%%%%%%%%%%%%%%%%%%%%%%%%%%%%%%%%%%%%%%%%%%%%%
%%%%%%%%%%%%%%%%%   subsection 5.4

{\bf 5.4. The case 2.2, $B$ is not contained in $Q_P$.} Now assume
that the case 2.2 takes place, where $B\not\subset Q_P$. Now, if
$B$ is not contained in a quadric, we obtain a contradiction,
arguing as in Subsection 5.1. If $B$ is contained in a quadric,
but not contained in a hyperplane, then we obtain a contradiction,
arguing as in Subsection 5.2. Therefore we assume that
$B\subset\Pi$, where $\Pi\subset E_P$ is some hyperplane. Now, if
$M\geq 18$ or if $B^*\not\subset Q_P$, then we obtain a
contradiction in word for word the same way as in Subsection 5.3.
Therefore we assume that $M\leq 17$ and $B^*\subset Q_P$ is a
subvariety of codimension 2.\vspace{0.1cm}

Let us consider the pencil of hyperplane sections $|H_P-\Pi|$. Its
base set $\Delta=\mathop{\rm Bs}|H_P-\Pi|$ is a hyperplane section
of the tangent section $T_P$. Write
$$
Z_P=a\Delta+Z_*,
$$
where $a\in{\mathbb Z}_+$ and $Z_*$ does not contain $\Delta$ as a
component. For the subvariety $\Delta$ we have:
$$
\mathop{\rm deg}\Delta=M,\,\, \mathop{\rm
mult}\nolimits_o\Delta=2, \,\, \mathop{\rm
mult}\nolimits_B\Delta^+=0
$$
and $\mathop{\rm mult}_{B^*}\Delta^+=1$. Therefore for the cycle
$Z_*$ we have: $\mathop{\rm deg}Z_*=(4n^2-a)M$,
$$
\mathop{\rm mult}\nolimits_oZ_*= \mathop{\rm
mult}\nolimits_oZ_P-2a>8n^2-2a, \quad \mathop{\rm
mult}\nolimits_BZ^+_*=\mathop{\rm mult}\nolimits_BZ^+_P
$$
and $\mathop{\rm mult}_{B^*}Z^+_*=\mathop{\rm
mult}_{B^*}Z^+_P-a>4n^2-a$. Now let $R\in|H_P-\Pi|$ be a general
divisor. By construction, $R$ does not contain irreducible
components of the cycle $Z_*$ and for that reason the effective
cycle $(Z_*\circ R)$ of codimension 2on $R$ is well defined. Let
$Z_R=(Z_*\circ R)=Z_0+Z_1$ be the $T_R$-decomposition of the cycle
$Z_R$.\vspace{0.1cm}

Let $Q_R=Q_P\cap\Pi=T^+_R\cap E_P$ be the (projectivized) tangent
cone to $T_R$, a quadric in $E_R=R^+\cap E_P=\Pi$. The subvariety
$B^*$ is a {\it prime divisor} on $Q_R$ and for that reason is cut
out on $Q_R$ by a hypersurface in $E_R$ of degree $\delta^*\geq
1$, so that $2\delta^*=d^*=\mathop{\rm deg}B^*$.\vspace{0.1cm}

{\bf Proposition 5.5.} {\it The equality $\delta^*=1$ holds, that
is, $B^*$ is a hyperplane section of $Q_R$.}\vspace{0.1cm}

{\bf Proof.} Assume the converse: $\delta^*\geq 2$. In that case
$d^*\geq 4$. Therefore, the inequality
\begin{equation}\label{07.10.2013.2}
\mathop{\rm mult}\nolimits_oZ_R\geq 4\mathop{\rm
mult}\nolimits_{B^*}Z^+_R
\end{equation}
holds. To compute the left hand part, write
$$
(Z^+_*\circ R^+)=Z^+_R+\beta B+N,
$$
where $N$ is an effective divisor on $\Pi$, not containing $B$ as
a component and $\beta\geq\mu_B$. By the intersection theory,
$$
\mathop{\rm mult}\nolimits_oZ_R= \mathop{\rm
mult}\nolimits_oZ_*+\beta d_B+d_N,
$$
where $d_N=\mathop{\rm deg}N$. On the other hand, by assumption
$B^*$ is not contained in a hyperplane $\Pi$, that is, $\langle
B^*\rangle=\Pi$ and for that reason the inequalities
$$
2\mathop{\rm mult}\nolimits_{B^*}N\leq d_N\quad \mbox{и}\quad
2\mathop{\rm mult}\nolimits_{B^*}B\leq d_B
$$
hold. Therefore, we have the estimate
$$
\mathop{\rm mult}\nolimits_{B^*}Z^+_R\geq \mathop{\rm
mult}\nolimits_{B^*}Z^+_*-\frac12\beta d_B-\frac12d_N.
$$
Besides, we remember that the inequalities
$$
\frac{\mathop{\rm mult}_o}{\mathop{\rm deg}}Z_0\leq
\frac{8}{3(M-2)}\quad \mbox{и} \quad \frac{\mathop{\rm
mult}_o}{\mathop{\rm deg}}Z_1\leq\frac{2}{M-2}
$$
hold, and also the inequalities $\mu>8n^2$ and $\mu+\mu_B>12n^2$.
Using MAPLE, it is easy to check (replacing, as usual, strict
inequalities by non-strict ones), that the system of linear
equations and inequalities, obtained above, has no solutions for
$M\geq 13$. Proof of proposition 5.5 is complete.\vspace{0.1cm}

Therefore, $B^*=\Theta\cap Q_R$, where $\Theta\subset\Pi=E_R$ is a
hyperplane. Instead of the inequality (\ref{07.10.2013.2}) we have
a weaker estimate
$$
\mathop{\rm mult}\nolimits_oZ_R\geq 2\mathop{\rm
mult}\nolimits_{B^*}Z^+_R
$$
and it is no longer sufficient to obtain a contradiction. Let us
consider the linear system $|H_R-\Theta|$ on $R$ and set
$\Delta^*=\mathop{\rm Bs}|H_R-\Theta|$ to be its base set (a
divisor on the tangent section). By the regularity conditions we
have $\mathop{\rm mult}_{B^*}(\Delta^*)^+=1$. Write down
$Z_0=c\Delta^*+Z_{\sharp}$, where $c\in{\mathbb Z}_+$ and
$Z_{\sharp}$ does not contain $\Delta^*$ as a component. For a
general divisor $D\in|H_R-\Theta|$ the effective cycle $(D\circ
Z_{\sharp})$ of codimension two on $T_R$ is well defined and
satisfies the inequalities
$$
\mathop{\rm mult}\nolimits_o(D\circ Z_{\sharp})\geq \mathop{\rm
mult}\nolimits_oZ_{\sharp}+ 2\mathop{\rm
mult}\nolimits_{B^*}Z^+_{\sharp}
$$
and
$$
\frac{\mathop{\rm mult}_o}{\mathop{\rm deg}}(D\circ Z_{\sharp})
\leq\frac{4}{M-2}.
$$
Adding the corresponding normalized inequalities to the previous
ones and using MAPLE, we see that for $M\geq 13$ the case under
consideration is impossible.\vspace{0.1cm}

This completes the exclusion of the case 2.2 under the assumption
that $B\not\subset Q_P$.\vspace{0.3cm}

%%%%%%%%%%%%%%%%%%%%%%%%%%%%%%%%%%%%%%%%%%%%%%%%%%%%%%%%%%%%%%%%%%
%%%%%%%%%%%%%%%%%%%%%%%%%%%%%%%%%%%%%%%%%%%%%%%%%%%%%%%%%%%%%%%%%%
%%%%%%%%%%%%%%%%%%%%   subsection 5.5

{\bf 5.5. The case 2.2, $B$ is contained in $Q_P$.} Assume that
$B\subset Q_P$. Note, first of all, that $B$ is not contained in a
hyperplane (that is, it is not a hyperplane section of $Q_P$):
such an option is excluded by word for word the same arguments as
those that were used in the case $B\not\subset Q_P$,
$B\subset\Pi$, where $\Pi\subset E_P$ is a hyperplane. In
particular, $d_B\geq 4$ and the estimate (\ref{12.09.1}) holds
.\vspace{0.1cm}

{\bf Proposition 5.6.} {\it The subvariety $B^*\subset Q_P$ of
codimension two is contained in a hyperplane} $\Pi\subset
E_P$.\vspace{0.1cm}

{\bf Proof.} Assume the converse. Let $\Lambda\subset Q_P$ be a
general linear subspace of maximal dimension,
$B^*_{\Lambda}=B^*\cap\Lambda\subset\Lambda$ an irreducible
subvariety of codimension two. For the linear span $\langle
B^*_{\Lambda}\rangle$ there are three options:\vspace{0.1cm}

1) $\langle B^*_{\Lambda}\rangle=\Lambda$,\vspace{0.1cm}

2) $\langle B^*_{\Lambda}\rangle$ is a hyperplane in
$\Lambda$,\vspace{0.1cm}

3) $\langle B^*_{\Lambda}\rangle=B^*_{\Lambda}$ is a subspace of
codimension two in $\Lambda$.\vspace{0.1cm}

Note at once, that the third option does not realize: 3) implies
that $\mathop{\rm deg}B^*=2$ and then $B^*$ is contained in a
hyperplane, contrary to our assumption.\vspace{0.1cm}

Furthermore, $\mathop{\rm Sec}(B^*_{\Lambda})=\langle
B^*_{\Lambda}\rangle$. Set
$$
W=\overline{\bigcup_{\Lambda\subset Q_P} \langle
B^*_{\Lambda}\rangle}.
$$
It follows from what was said that either $W$ is an irreducible
divisor on $Q_P$, or $W=Q_P$. However, in the first case
$W\cap\Lambda$ is an irreducible hypersurface in $\Lambda$ (for a
general $\Lambda$) and for that reason $W\cap\Lambda=\langle
B^*_{\Lambda}\rangle$ is a hyperplane in $\Lambda$, and then $W$
is a hyperplane section of the quadric $Q_P$, where $B^*\subset
W$, contrary to our assumption. Therefore, $W=Q_P$. From here we
get the following fact.\vspace{0.1cm}

{\bf Lemma 5.3.} {\it For any effective divisor $Y$ on the quadric
$Q_P$ the inequality}
$$
\mathop{\rm deg}Y\geq 4\mathop{\rm mult}\nolimits_{B^*}Y
$$
{\it holds (the degree $\mathop{\rm deg}Y$ is understood as the
degree of an effective cycle of codimension 2 on
$E_P$).}\vspace{0.1cm}

{\bf Proof.} Denote by the symbol $H_Q$ the class of a hyperplane
section of the quadric $Q_P$, so that $Y\sim \gamma H_Q$ for some
$\gamma\geq 1$, where $\mathop{\rm deg}Y=2\gamma$. Let $\Lambda$
be a general linear subspace of maximal dimension on $Q_P$ and
$L\subset\Lambda$ a general secant line of the variety
$B^*_{\Lambda}$. Since the lines $L$ sweep out $Q_P$, we may
assume that $L\not\subset|Y|$. Let $x,y\in B^*_{\Lambda}$ be
general points, where $L=[x,y]$. We have
$$
(L\cdot Y)_{Q_P}= \gamma\geq (L\cdot Y)_x+(L\cdot Y)_y\geq
2\mathop{\rm mult}\nolimits_{B^*}Y,
$$
when the claim of the lemma follows. Q.E.D.\vspace{0.1cm}

Now let $Z_P=Z_0+Z_1$ be, as usual, the $T_P$-decomposition of the
cycle $Z_P$. Setting $\lambda_i=\frac{1}{n^2}\mathop{\rm
mult}\nolimits_BZ_i$, $i=0,1$, we obtain the following system of
linear equations and inequalities:
(\ref{25.06.1},\ref{25.06.2},\ref{25.06.4}), and also the estimate
\begin{equation}\label{07.10.2013.1}
\mu_0+\mu_1+\lambda_0+\lambda_1>10+2\sqrt{2}
\end{equation}
instead of (\ref{25.06.3}), and also the estimates
\begin{equation}\label{23.10.2013.2}
2\mu_1+d_B\lambda_1\leq\frac{4M}{M-1}d_1,
\end{equation}
\begin{equation}\label{23.10.2013.1}
\mu_0\geq d_B\lambda_0,\quad \mu_1\geq d_B\lambda_1.
\end{equation}
Now set $\xi_i=\frac{1}{n^2}\mathop{\rm
mult}\nolimits_{B^*}Z^+_i$, $i=0,1$. By the lemma shown above, the
estimate
$$
\mu_0\geq 4\xi_0
$$
holds, besides, $\mu_1\geq\xi_1$ and, as we know, $\xi_0+\xi_1>4$.
The inequality (\ref{23.10.2013.2}) can be sharpened. Write down
$$
(Z^+_1\circ T_P)=(Z_1\circ T_P)^++N,
$$
where $N$ is an effective divisor on the quadric $Q_P$. Set
$d_N=\frac{1}{n^2}\mathop{\rm deg}N$, then we get
$$
d_N\geq d_B\lambda_1
$$
and the estimate
$$
2\mu_1+d_N\leq\frac{4M}{M-1}d_1
$$
holds. Setting $\xi_N=\frac{1}{n^2}\mathop{\rm mult}_{B^*}N$ and
applying Lemma 5.3, we obtain the inequality
\begin{equation}\label{20.09.1}
d_N\geq 4\xi_N.
\end{equation}
Obviously,
$$
\frac{1}{n^2}\mathop{\rm mult}\nolimits_{B^*}(Z_1\circ T_P)^+
\geq\xi_1-\xi_N,
$$
so that, applying Lemma 5.3 once again, we get the inequality
$$
2\mu_1+d_N\geq 4(\xi_1-\xi_N).
$$
Using MAPLE, we check that the system of linear equations and
inequalities, obtained above, is incompatible. Q.E.D. for
Proposition 5.6.\vspace{0.3cm}

%%%%%%%%%%%%%%%%%%%%%%%%%%%%%%%%%%%%%%%%%%%%%%%%%%%%%%%%%%%%%%%%
%%%%%%%%%%%%%%%%%%%%%%%%%%%%%%%%%%%%%%%%%%%%%%%%%%%%%%%%%%%%%%%%
%%%%%%%%%%%%%%%%%%   subsection 5.6

{\bf 5.6. Exclusion of the case 2.2.} Now let us assume that the
hyperplane $\Pi\supset B^*$ is the only hyperplane in $E_P$ with
that property, that is, $B^*$ is not the intersection of $Q_P$
with a linear subspace $\Theta\subset E_P$ of codimension two. In
particular, $d^*=\mathop{\rm deg}B^*\geq 4$. Let $R\in|H_P-\Pi|$
be a general divisor of the pencil. Write down $Z_P=a\Delta+Z_*$,
where $\Delta=\mathop{\rm Bs}|H_P-\Pi|$, $a\in{\mathbb Z}_+$ and
$Z_*$ does not contain $\Delta$ as a component. To simplify the
formulas, we will assume that $a=0$ and $Z_*=Z_P$: if $a\geq 1$,
then the system of linear equations and inequalities, obtained
below, remains incompatible, which is easy to check.\vspace{0.1cm}

So $Z_P=Z_0+Z_1$ is the $T_P$-decomposition of the cycle $Z_P$ and
$\Delta$ is not an irreducible component of the cycle $Z_0$.
Setting, as usual,
$$
\mu_i=\frac{1}{n^2}\mathop{\rm mult}\nolimits_oZ_i,\quad
\lambda_i=\frac{1}{n^2} \mathop{\rm mult}\nolimits_BZ^+_i,
$$
and $d_i=\frac{1}{Mn^2}\mathop{\rm deg}Z_i$, $i=0,1$, we obtain
the standard set of linear equations and inequalities:
(\ref{25.06.1},\ref{25.06.2},\ref{07.10.2013.1}), and also the
inequalities (\ref{25.06.4}), (\ref{23.10.2013.1}) with $d_B=4$
and the estimate
$$
\mu_1\leq\frac{2M}{M-1}d_1.
$$
Set $\xi_i=\frac{1}{n^2}\mathop{\rm mult}_{B^*}Z^+_i$, $i=0,1$. In
our case $\xi_0+\xi_1>4$.\vspace{0.1cm}

{\bf Lemma 5.4.} {\it The following inequality holds:}
$$
\mu_0\geq 2(\lambda_0+\xi_0).
$$

{\bf Proof.} $Z^+_0$ is an effective divisor on $T^+_P$, and its
projectivized tangent cone $Z^+_0\cap E_P$ is an effective divisor
on the quadric $Q_P$. Let $\Lambda\subset Q_P$ be a general linear
subspace. Let $p\in B^*\cap\Lambda$ and $q\in B\cap\Lambda$ be
points of general position. The lines $L=[pq]$ sweep out $\Lambda$
and for that reason we may assume that $L\not\subset Z^+_0$.
Therefore, for the intersection numbers on $Q_P$ we have:
$$
\frac12\mu_0n^2=(L\cdot(Z^+_0\cap E_P))_{Q_P}\geq(L\cdot(Z^+_0\cap
E_P))_p+ (L\cdot(Z^+_0\cap E_P))_q\geq (\xi_0+\lambda_0)n^2,
$$
which is what we claimed. Q.E.D. for the lemma.\vspace{0.1cm}

Now write down
$$
(Z^+_i\circ R^+)=(Z_i\circ R)^++N_i,
$$
where $N_1$ is an effective divisor on $\Pi$, and $N_0=c(\Pi\cap
Q_P)$, $c\in{\mathbb Z}_+$. Since $\langle B^*\rangle=\Pi$, the
inequality
$$
\mathop{\rm deg}N_i\geq 2\mathop{\rm mult}\nolimits_{B^*}N_i
$$
holds (for $i=0$ it is the equality, since obviously $\mathop{\rm
mult}_{B^*}N_0=1$). Setting
$$
\zeta_i=\frac{1}{n^2}\mathop{\rm mult}\nolimits_{B^*}(Z_i\circ
R)^+,
$$
we obtain inequalities
$$
\zeta_i\geq\xi_i-\frac12n_i,
$$
where $n_i=\frac{1}{n^2}\mathop{\rm deg}N_i$. Setting
$\alpha_i=\frac{1}{n^2}\mathop{\rm mult}_o(Z_i\circ R)$, $i=0,1$,
we obtain the set of standard estimates
$$
\alpha_i\geq\mu_i+n_i,\quad
\alpha_0\leq\mathop{\rm max}\left(3,\frac{8M}{3(M-2)}\right)\cdot d_0,\quad
\alpha_i\geq 4\zeta_i,\quad i=0,1
$$
(the last is true by the inequality $\mathop{\rm deg}B^*\geq 4$,
as $\langle B^*\rangle=\Pi$). Besides, one more important
inequality holds.\vspace{0.1cm}

{\bf Lemma 5.5.} {\it The following estimate holds:}
$$
4n_0\geq 4\xi_0-\mu_0.
$$

{\bf Proof.} Once again, let $\Lambda\subset Q_P$ be a general
linear subspace of the maximal dimension and $L$ a general secant
line of the variety $B^*\cap\Lambda$. The lines $L$ sweep out the
hyperplane section $\Pi\cap Q_P$ and for that reason it is
sufficient to show the inequality
\begin{equation}\label{21.09.1}
\beta=\frac{1}{n^2}\mathop{\rm mult}\nolimits_LZ^+_0\geq
\frac14(2\xi_0-\frac12\mu_0).
\end{equation}
Since $n_0\geq 2\beta$, the inequality (\ref{21.09.1}) implies the
claim of our lemma.\vspace{0.1cm}

Consider a general line $L^*\subset\Lambda$, intersecting $L$, and
let $S\ni o$ be a generic two-dimensional germ of an isolated
quadratic singularity at the point $o$, $S\subset T_P$, such that
$S^+\cap E_P=L+L^*$, and $S^+$ is a non-singular surface.
Obviously,
$$
Z^+_0|_{S^+}\sim -\left(\frac12\mathop{\rm
mult}\nolimits_oZ_0\right)E_P|_{S^+},
$$
whereas the effective 1-cycle $Z^+_0|_{S^+}$ has the line $L$ as a
component of the multiplicity $\beta n^2$. Taking this component
out, we obtain that the effective 1-cycle
$$
C=(Z^+_0|_{S^+}-\beta n^2L)
$$
does not have $L$ as a component, and its multiplicity at two
distinct points $p,q\in L$ is at least
$$
(\xi_0-\beta)n^2.
$$
Computing the intersection $(C\cdot L)$, we obtain the inequality
$$
\frac12\mu_0+2\beta\geq 2(\xi_0-\beta),
$$
which is what we need. Q.E.D. for the lemma.\vspace{0.1cm}

Finally, adding the inequality of Lemma 5.5 to the previous
estimates, we obtain an incompatible system of linear equations
and inequalities (checked using MAPLE), which completes the
exclusion of the case under consideration.\vspace{0.1cm}

Therefore, the only remaining possibility is when $B^*=\Theta\cap
Q_P$, where $\Theta\subset E_P$ is a linear subspace of
codimension two. The claim of Lemma 5.4 is valid. Let
$R\in|H_P-\Theta|$ be a general divisor, $Z_R=(Z_P\circ R)$ an
effective cycle of codimension two on $R$, $Z_P=Z_0+Z_1$ is, as
usual, the $T_P$-decomposition of the cycle $Z_P$. We get the
standard set of linear equalities and inequalities for that
decomposition:
(\ref{25.06.1},\ref{25.06.2},\ref{25.06.4},\ref{07.10.2013.1},
\ref{23.10.2013.1}) and (\ref{23.10.2013.2}) with $d_B=4$. Now let
us consider the cycle $Z_R$ more carefully. Set
$$
(Z_i\circ R)=Z^{\sharp}_i+c_in^2\Delta,
$$
where $\Delta=\mathop{\rm Bs}|H_P-\Theta|$ is a hyperplane section
of the hypersurface $T_R$; note that none of the components of the
cycle $Z^{\sharp}_1$ is not contained in $T_R$. The support of the
cycle $Z^{\sharp}_2$ is contained in $T_R$. For the subvariety
$\Delta$ we obviously have: $\mathop{\rm deg}\Delta=M$,
$\mathop{\rm mult}\nolimits_o\Delta=2$ and $\mathop{\rm
mult}\nolimits_{B^*}\Delta^+=1$. Obviously,
$$
\mathop{\rm mult}\nolimits_o(Z_1\circ R)=
\mathop{\rm mult}\nolimits_oZ_1,
\quad \mathop{\rm mult}\nolimits_{B^*}(Z_1\circ R)^+=
\mathop{\rm mult}\nolimits_{B^*}Z^+_1.
$$
Setting $\mu^{\sharp}_i=\frac{1}{n^2}\mathop{\rm
mult}_oZ^{\sharp}_i$ for $i=0,1$, we obtain the inequality
$$
\mu^{\sharp}_1+(\xi_1-c_1)\leq\frac{2M}{M-2}(d_1-c_1).
$$
Furthermore, the following equalities
$$
\mathop{\rm mult}\nolimits_o(Z_0\circ R)=\mathop{\rm
mult}\nolimits_oZ_0,\quad \mathop{\rm
mult}\nolimits_{B^*}(Z_0\circ R)^+= \mathop{\rm
mult}\nolimits_{B^*}Z^+_0
$$
hold. The cycle $Z^{\sharp}_0$ is an effective divisor on $T_R$,
which does not contain $\Delta$ as a component. Let
$R^*\in|H_P-\Theta|$ be another general divisor. Obviously,
$$
R^*\cap T_R=R\cap R^*\cap T_P=\Delta,
$$
so that none of the components of the cycle $Z^{\sharp}_0$ is not
contained in $R^*$ and therefore the cycle
$Z^*_0=(Z^{\sharp}_0\circ R^*)$ of codimension two on $T_R$ is
well defined. The cycle $Z^*_0$ is an effective divisor on
$\Delta$. Setting $\mu^*_0=\mathop{\rm mult}_oZ^*_0$, we obtain
the inequality
$$
\mu^*_0\geq\mu^{\sharp}_0+2(\xi_0-c_0).
$$
By the regularity conditions on the hypersurface $R$ the
inequality
$$
\mu^*_0\leq\frac{4M}{M-2}(d_0-c_0)
$$
holds. But it is not hard to obtain a stronger estimate. By the
regularity conditions and the Lefschetz theorem we have:
$$
(\Delta\circ T_o(T_R))=\Delta\cap T_o(T_R)=\Delta\cap T_o(T_P)
$$
is an irreducible reduced divisor on $\Delta$, which has the
degree $2M$ and the multiplicity precisely 6 at the point $o$. Let
$Y$ be an irreducible component of the cycle $Z^*_0$. Then either
the inequality $(\mathop{\rm mult}_o/\mathop{\rm deg})Y\leq(3/M)$
holds, or $Y$ is not contained in the hypertangent divisor
$T_o(T_P)$, so that the estimate
$$
\frac{\mathop{\rm mult}_o}{\mathop{\rm deg}}Y\leq\frac{10}{3(M-2)}
$$
is true. From this it follows, that the inequality
$$
\mu^*_0\leq\mathop{\rm
max}\left(3,\frac{10M}{3(M-2)}\right)(d_0-c_0)
$$
holds.\vspace{0.1cm}

Using MAPLE, it is easy to check that the system of linear
equations and inequalities for $\mu_*$, $d_*$, $c_*$,
$\mu^{\sharp}_*$ and $\mu^*_0$, obtained above, has no
solutions.\vspace{0.1cm}

The case 2.2 is completely excluded.\vspace{0.3cm}

%%%%%%%%%%%%%%%%%%%%%%%%%%%%%%%%%%%%%%%%%%%%%%%%%%%%%%%%%%%%%%%%%%%%%
%%%%%%%%%%%%%%%%%%%%%%%%%%%%%%%%%%%%%%%%%%%%%%%%%%%%%%%%%%%%%%%%%%%%%
%%%%%%%%%%%%%%   subsection 5.7

{\bf 5.7. Exclusion of the case 2.3.} Assume that the case 2.3
takes place. We have $\mu_B=\mathop{\rm mult}_BZ^+_P>4n^2$, so
that we get the following sequence of inequalities:
$$
12n^2\geq\mathop{\rm mult}\nolimits_oZ_P\geq d_B\mu_B>4d_Bn^2,
$$
where $d_B=\mathop{\rm deg}B\geq 2$ (the case of a linear subspace
was excluded by Proposition 4.1), whence we conclude that $d_B=2$,
that is, $B$ is a quadric in some hyperplane $\Pi\subset
E_P$.\vspace{0.1cm}

If $B\not\subset Q_P$, then we argue as in Subsection 5.3: we
write down $Z_P=Z_0+Z_1$ and intersect $Z_1$ with a general
divisor $R\in|H_P-\Pi|$. Since $\mu_B>4n^2$, we obtain the linear
inequalities
$$
\begin{array}{l}
\displaystyle\mathop{\rm
mult}\nolimits_oZ_0\leq\frac{3}{M}\mathop{\rm deg}Z_0,\\
\\
\displaystyle \mathop{\rm mult}\nolimits_oZ_1+8n^2<
\frac{8}{3(M-2)}\mathop{\rm deg}Z_1\leq\frac{4}{M}\mathop{\rm deg}Z_1\\
\end{array}
$$
which hold for $M\geq 6$. Putting together and recalling that
$\mathop{\rm mult}_oZ_P>8n^2$, we obtain a contradiction,
excluding the possibility $B\not\subset Q_P$.\vspace{0.1cm}

So let us assume that $B=\Pi\cap Q_P$ is a hyperplane section of
the quadric $Q_P$. Consider $\Delta=\mathop{\rm Bs}|H_P-\Pi|$,
which is a hyperplane section of the variety $T_P$. Write down
$$
Z_P=a\Delta+Z_*,
$$
where $a\in{\mathbb Z}_+$ and $Z_*$ does not contain $\Delta$ as a
component. For $\Delta$ we have: $\mathop{\rm deg}\Delta=M$,
$\mathop{\rm mult}_o\Delta=2$ and $\mathop{\rm mult}_B\Delta^+=1$.
Therefore, $\mathop{\rm deg}Z_*=(4n^2-a)M$, and for the
multiplicities we have the equalities
$$
\mathop{\rm mult}\nolimits_oZ_*=\mu-2a,\quad \mathop{\rm
mult}\nolimits_BZ^+_*=\mu_B-a.
$$
Now, arguing in word for word the same way as above in the case
$B\not\subset Q_P$, where $Z_P$ is replaced by $Z_*$, we obtain a
contradiction. The case 2.3 is excluded.\vspace{0.1cm}

Proof of Theorem 5 is complete.

\begin{flushleft}
Department of Mathematical Sciences,\\
The University of Liverpool
\end{flushleft}

\noindent{\it pukh@liv.ac.uk}

\end{document}